\documentclass[12pt]{article}
\usepackage{latexsym, amssymb, enumerate, amsmath}

\usepackage{graphicx}  
\newtheorem{thm}{Theorem}[section]
\newtheorem{prop}[thm]{Proposition}
\newtheorem{lem}[thm]{Lemma}

\newtheorem{defn}[thm]{Definition}

\begin{document}
\begin{center}
{\Large {\bf Bifurcations of 
Wavefronts on an  $r$-corner} }
\vspace*{0.4cm}\\
{\large Takaharu Tsukada}\footnote{Higashijujo 3-1-16 
Kita-ku, Tokyo 114-0001
JAPAN. e-mail : tsukada@math.chs.nihon-u.ac.jp}
\vspace*{0.2cm}\\
{\large  College of Humanities \& Sciences, Department of Mathematics,\\
 Nihon University}
\end{center}
\begin{abstract}
We introduce the notion of {\em reticular Legendrian
unfoldings}  in order to investigate stabilities and a genericity of bifurcations 
of wavefronts 
generated by a hypersurface germ with a boundary, a corner, or an 
$r$-corner in a smooth $n$ dimensional manifold.
We define several stabilities of reticular Legendrian
unfoldings and prove that they and the stabilities of corresponding 
generating families are all equivalent 
 and give the classification
of all generic 
bifurcations of wavefronts in the cases $r=0,n\leq 5$ and $r=1,n\leq 3$ 
respectively.
\end{abstract}

%

\section{Introduction}\label{intro}
\quad
In \cite{janich1} K.J\"anich explained the wavefront propagation mechanism 
on a manifold which is completely described by a positive and positively 
homogeneous {\em Hamiltonian function} on the cotangent bundle and investigated  
the local gradient models given by the ray length function. 
Caustics and Wavefronts generated by an initial wavefront 
which is a hypersurface germ
without boundary in the manifold were investigated as 
Lagrangian and Legendrian singularities 
by V.I.Arnold (cf., \cite{arnold:text}). 
%

In this paper and its
prequel \cite{tPKfunct},
we investigate the stabilities and the genericity of bifurcations of 
wavefronts generated by a hypersurface germ with {\em an $r$-corner}.
Wavefronts generated by all edges of the hypersurface at a time give {\em a contact regular $r$-cubic configuration}
on the $1$-jet bundle.
All wavefronts around a time give a one-parameter family of contact 
regular $r$-cubic configurations on the $1$-jet bundle.
In order to consider such families,
we shall introduce the notion of
{\em unfolded contact regular $r$-cubic configurations} on the big 
$1$-jet bundle.
A wavefront of an unfolded contact regular $r$-cubic configuration
is the big front of the corresponding one-parameter family of contact regular $r$-cubic configurations. 
We shall consider their generating families and equivalence relations. 
\begin{figure}[htbp]
 \begin{minipage}{0.30\hsize}
  \begin{center}
    \includegraphics[width=3cm,height=3cm]{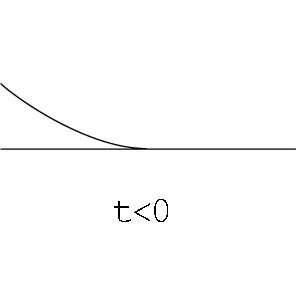}
  \end{center}
 \end{minipage}
$\leftrightarrow $
 \begin{minipage}{0.30\hsize}
   \begin{center}
     \includegraphics[width=3cm,height=3cm]{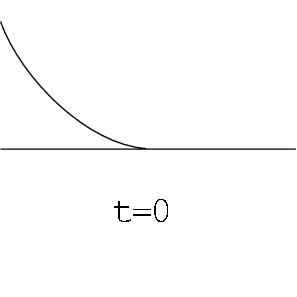}
   \end{center}
 \end{minipage}
$\leftrightarrow $
 \begin{minipage}{0.30\hsize}
  \begin{center}
  \includegraphics[width=3cm,height=3cm]{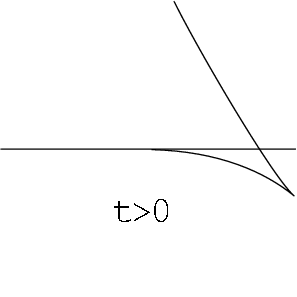}
  \end{center}
\end{minipage}

 \caption{The generic bifurcation  of wavefronts on a boundary ${}^1B_3$}
\end{figure}

In order to investigate the stabilities and the genericity of
unfolded contact regular $r$-cubic configurations,
we introduce the notion of {\em reticular 
Legendrian unfoldings} which is
a generalized notion of {\em Legendrian unfoldings} given 
by  S.Izumiya (cf., \cite{izumiya1}) for our situation.
We shall define several stabilities of reticular 
Legendrian unfoldings
and prove that they and the stabilities of their generating families are all 
equivalent.
We shall also classify generic reticular Legendrian unfoldings and
give all figures of their wavefronts in the case $r=1,n\leq 3$.

In history, our theory in the case $r=0$ is investigated as the theory
of Legendrian unfoldings by S.Izumiya (cf., \cite{izumiya1}) and
the classification list of generic Legendrian unfoldings is given by using
V.M.Zakalyukin's theory (cf., \cite{zaka1}) in which he has classified
generic {\em quasihomogeneous function germs} under 
{\em the $t$-${\cal P}$-${\cal K}$-equivalence}.
In \cite{tPKfunct} we have classified not only quasihomogeneous function germs 
but also all smooth function germs under {\em the reticular $t$-${\cal P}$-${\cal K}$-equivalence}
which is a generalized relation of  the $t$-${\cal P}$-${\cal K}$-equivalence.

I.G.Scherbak has studied the theory of boundary fronts in \cite{IScher}
and this corresponds to our theory in the case $r=1$.
She has introduced the notion of Legendrian pairs which is 
corresponding to the notion of regular $1$-cubic configurations.
But they are not strictly formulated and no proof is given. 
In this paper, we shall define and prove our theory strictly.
Since her equivalence relation and our one of function germs
slightly different to each other (see the remark in Section \ref{unfold:sec}),
the figures of wavefronts ${}^0B_2$ of fig.2 and ${}^1B_3$ of fig.3(cf., Figure 1), and 
${}^0B_3$ of fig.4 in \cite[p.365]{IScher} do not coincide with our figures.
The classification list of function germs also different to each other
(compare \cite[p.371 Proposition 3]{IScher} with Theorem \ref{classfunct:tth}).

This paper consists of six sections. In Section \ref{pre:sec} we shall give 
the motivation for this paper and give the review of
stabilities under the  reticular $t$-${\cal P}$-${\cal K}$-equivalence relation
of function germs which play important roles as 
generating families of reticular Legendrian unfoldings (cf., \cite{tPKfunct}).
We shall also give the review of the theory of contact regular $r$-cubic configurations.
In Section \ref{RLeg:unfo} we shall introduce the notion of 
reticular Legendrian unfoldings and consider their generating families.
In Section \ref{map:sec} we shall investigate several stabilities of 
reticular Legendrian unfoldings.
In Section \ref{generic:sec} we shall reduce our investigation to 
finite dimensional jet spaces and give the classification of generic 
reticular Legendrian unfoldings in the cases $r=0$ and $1$ respectively.
In Section \ref{classover1} we shall show that our method do not work well
for the cases $r\geq 2$ because of modalities of generating families.
All maps considered here are differentiable of class $C^\infty$ unless stated
otherwise.
\section{Preliminary}\label{pre:sec}
\subsection{Propagation mechanism of wavefronts}
\quad 
Let us start with the propagation mechanism of wavefronts generated by a hypersurface germ $V^0$ 
with an $r$-corner in an $(n+1)(=r+k+1)$-dimensional smooth manifold $M$ which is given in \cite{janich1}.
Let $H:T^*M \backslash 0 \rightarrow {\mathbb R}$ be a fixed Hamiltonian function, which we suppose that $H(\lambda \xi)=\lambda H(\xi)$ for all $\lambda >0$ and
 $\xi \in T^*M \backslash 0$.
For example, consider a Riemann manifold $M$ and $H$ be length of covectors
in $T^*M \backslash 0$.

The manifold $E=H^{-1}(1)$ has the contact structure defined by
the restriction of the canonical $1$-form on $T^*M$ and 
the projection $\pi:E\rightarrow M$ is a spherical cotangent bundle.

Let ${\mathbb H}^r=\{ (x_1,\cdots,x_r)\in {\mathbb R}^r|x_1\geq 0,\cdots,x_r\geq 0\}$ 
be an $r$-corner, $\xi_0\in E$, $t_0\geq 0$, and $V^0$ be the initial hypersurface germ 
defined by the image of the immersion 
$\iota:({\mathbb H}^r\times {\mathbb R}^k,0)\rightarrow M$ such that 
$\iota(0)=\pi(\xi_0),\ \xi_0|_{T_{\iota(0)}V^0}=0$.
Let $\eta_0$ be the image of the phase flow of the Hamiltonian vector field
$X_H$ at $(t_0,\xi_0)$.
Since the flow preserves values of $H$ and the contact structure on $E$, it induces the
contact embedding germs
$C_t:(E,\xi_0)\rightarrow E$ for $t$ around $t_0$ which depends smoothly on $t$.
We define the $\sigma$-edge $V^0_\sigma$ of $V^0$ by 
$V^0\cap \{ x_\sigma =0 \}$ for $\sigma \subset I_r=\{1,\ldots,r\}$.
Let $L^0_\sigma$ be the initial covectors in $E$ generated by $V^0_\sigma$
to conormal directions, that is
\[ L^0_\sigma=\{ \xi_q\in E \ | \  q\in V^0_\sigma,\ \xi_q|_{T_qV^0_\sigma}=0 \}. \]
We may regard $L^0_\sigma$ as the lift of $V^0_\sigma$.
We also define that $L_{\sigma,t}=C_t(L^0_\sigma)$ for $\sigma\subset I_r,
t\in ({\mathbb R},t_0)$.
Then the wavefront $W_{\sigma,t}$ generated by $V^0_\sigma$
to conormal directions at time $t$ is given by  
$W_{\sigma,t}=\pi(L_{\sigma,t})$ for $\sigma \subset I_r, t\in ({\mathbb R},t_0)$.
\begin{figure}[htbp]
  \begin{center}
    \includegraphics[width=16cm,height=4cm]{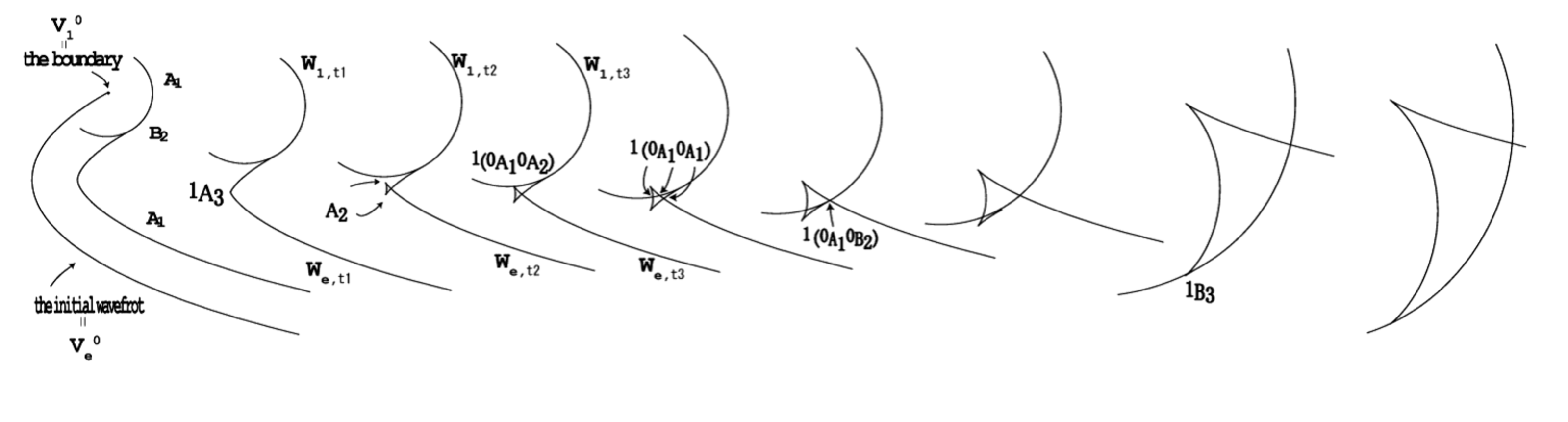}
  \end{center}
\caption{The initial wavefront $V^0$ with the boundary and 
the generated wavefronts {\rm(} $e=\emptyset$, $t_1<t_2<t_3${\rm)}}\label{pre:fig}
\end{figure}  

We are concerned with the stabilities and the genericity of 
bifurcations of wavefronts 
$\{ W_{\sigma,t}\}_{\sigma\subset I_r}$ for $t$ around $t_0$
with respect to perturbations of $V^0$.

Since we shall discuss local situations,  
we may  identify $(E,\xi_0)$ with $(J^1({\mathbb R}^{n},
{\mathbb R}),0)$, $\pi:(E,\eta_0)\rightarrow (M,\pi(\eta_0))$ with 
$\pi:(J^1({\mathbb R}^{n},{\mathbb R}),0)\rightarrow ({\mathbb R}^{n}\times {\mathbb R},0)$, and 
$t_0=0$,
where $\pi:J^1({\mathbb R}^{n},{\mathbb R})\rightarrow 
{\mathbb R}^{n}\times {\mathbb R}$ is the natural Legendrian bundle which 
is introduced in Section \ref{conta:sec}.
Then $C_t$ is identified 
with $C_t:(J^1({\mathbb R}^{n},{\mathbb R}),0)\rightarrow
J^1({\mathbb R}^{n},{\mathbb R})$ for $t\in ({\mathbb R},0)$ with $C_0(0)=0$.
The contact coordinate system on $(E,\xi_0)$ may be chosen that
$L^0_\sigma$ is given:
\[ L^0_\sigma=
\{ (q,z,p) \in (J^1({\mathbb R}^n,{\mathbb R}),0)|q_\sigma=p_{I_r-\sigma}=
q_{r+1}=\cdots=q_n=z=0,q_{I_r-\sigma}\geq 0 \} \mbox{ for each } \sigma \subset I_r.\]

We shall consider one-parameter families of contact regular $r$-cubic configurations 
$\{ L_{\sigma,t}\}_{\sigma\subset I_r,t\in ({\mathbb R},0)}$. 

\subsection{Stabilities of unfoldings}\label{unfold:sec}
\quad 
We review the main results of 
the theory of function germs with respect to 
the reticular $t$-${\cal P}$-${\cal K}$-equivalence relation
given in \cite{tPKfunct}.

We denote by ${\cal E}(r;k_1,s;k_2)$ the set of all germs at $0$ in 
${\mathbb H}^r\times {\mathbb R}^{k_1}$ of
smooth maps ${\mathbb H}^r\times {\mathbb R}^{k_1} \rightarrow 
{\mathbb H}^s\times {\mathbb R}^{k_2}$ and set ${\mathfrak M}(r;k_1,s;k_2)=
\{ f\in {\cal
E}(r;k_1,s;k_2)|f(0)=0 \}$.
We denote ${\cal E}(r;k_1,k_2)$ for ${\cal E}(r;k_1,0;k_2)$ and 
denote ${\mathfrak M}(r;k_1,k_2)$ for ${\mathfrak M}(r;k_1,0;k_2)$.

 If $k_2=1$ we write simply ${\cal E}(r;k)$ for 
${\cal E}(r;k,1)$
and ${\mathfrak M}(r;k)$ for ${\mathfrak M}(r;k,1)$. 
Then ${\cal E}(r;k)$ is an ${\mathbb R}$-algebra in the usual
way and ${\mathfrak M}(r;k)$ is its unique maximal ideal. 
We also denote by ${\cal E}(k)$ for 
${\cal E}(0;k)$
and ${\mathfrak M}(k)$ for ${\mathfrak M}(0;k)$.

We
denote by $J^l(r+k,p)$ the set of $l$-jets at $0$ of germs in ${\cal
E}(r;k,p)$. There are natural projections:
\[ \pi_l:{\cal E}(r;k,p)\longrightarrow J^l(r+k,p),
\pi^{l_1}_{l_2}:J^{l_1}(r+k,p)\longrightarrow J^{l_2}(r+k,p)\ (l_1 > l_2).  \]
We write $j^lf(0)$ for $\pi_l(f)$ for each $f\in {\cal E }(r;k,p)$.

Let $(x,y)=(x_1,\cdots,x_r,y_1,\cdots,y_k)$ be a fixed
coordinate system of $({\mathbb H}^r\times {\mathbb R}^k,0)$. 
We denote by 
${\cal B}(r;k)$ the group of
diffeomorphism germs $({\mathbb H}^r\times {\mathbb R}^{k},0)\rightarrow 
({\mathbb H}^r\times {\mathbb R}^{k},0)$ of the form:
\[ \phi(x,y)=(x_1\phi_1^1(x,y),\cdots,x_r\phi_1^r(x,y),\phi_2^1(x,y),\cdots,\phi_2^k(x,y)
). \]

We say that $f_0,g_0\in{\cal E}(r;k)$ are {\em reticular ${\cal K}$-equivalent} if
there exist $\phi\in{\cal B}(r;k)$ and a unit $a\in {\cal E}(r;k)$
such that $g_0=a\cdot f_0\circ \phi$.
We call $(\phi,a)$ {\em a reticular ${\cal K}$-isomorphism from $f_0$ to $g_0$}. \\

\noindent {\bf Remark}: The corresponding equivalence relation of function germs is 
given by I.G.Scherbak(see \cite[p.366 \S 2]{IScher}) as follows:
Function germs $f_0,g_0\in{\cal E}(1+k)$ are equivalent if there exist a 
diffeomorphism germ $\phi$ on $({\mathbb R}^{1+k},0)$ of the form
$\phi(x,y)=(x\phi_1(x,y),\phi_2^1(x,y),\ldots,\phi_2^k(x,y))$ and a unit
$a\in {\cal E}(r;k)$
such that $g_0=a\cdot f_0\circ \phi$.
The variable $x$ is defined on $({\mathbb H},0)$ 
and $\phi_1(0)>0$ in
our equivalence relation. 
On the other hand $x$ is defined on $({\mathbb R},0)$ and 
the condition $\phi_1(0)>0$ is not required in I.G.Scherbak's
 equivalence relation. 
These differences appear in the definition of wavefronts and consequently 
the figures of wavefronts in \cite{IScher} and the figures in this paper are different to each other.
\\

We denote by 
${\cal B}_n(r;k+n)$ the group of
diffeomorphism germs $({\mathbb H}^r\times {\mathbb R}^{k+n},0)\rightarrow 
({\mathbb H}^r\times {\mathbb R}^{k+n},0)$ of the form:
\[ \phi(x,y,u)=(x_1\phi_1^1(x,y,u),\cdots,x_r\phi_1^r(x,y,u),\phi_2^1(x,y,u),\cdots,\phi_2^k(x,y,u)
,\phi_3^1(u),\ldots,\phi_3^n(u)).\]
We denote $\phi(x,y,u)=(x\phi_1(x,y,u),\phi_2(x,y,u),\phi_3(u))$ and 
denote other notations analogously.\\

We say that $f,g\in{\cal E}(r;k+n)$ are {\em reticular ${\cal P}$-${\cal K}$-equivalent} if
there exist $\Phi\in{\cal B}_n(r;k+n)$ and a unit $a\in {\cal E}(r;k+n)$
such that $g=a\cdot f\circ \Phi$.
We call $(\Phi,a)$ {\em a reticular ${\cal P}$-${\cal K}$-isomorphism from $f$ to $g$}. \\

In convenience, we denote an unfolding of  a function germ 
$f(x,y,u)\in {\mathfrak M}(r;k+n)$ by
$F(x,y,t,u)\in {\mathfrak M}(r;k+m+n)$.

We
say that $F(x,y,t,u),G(x,y,t,u)\in {\cal E}(r;k+m+n)$
 are {\em reticular $t$-${\cal P}$-${\cal K}$-equivalent }
if there exist $\Phi$ of ${\cal B}(r;k+m+n)$ and a unit $a
\in{\cal E}(r;k+m+n)$ such that\\ 
(1) $\Phi$ can be written in the form:
\begin{equation}
\Phi(x,y,t,u)=(x_1\phi_1(x,y,t,u),\phi_2(x,y,t,u),\phi_3(t),\phi_4(t,u)),\label{pkreqn}
\end{equation}
(2) $G(x,y,t,u)=a(x,y,t,u) \cdot F\circ \Phi(x,y,t,u)$ for all 
$(x,y,t,u)\in ({\mathbb H}^r\times{\mathbb R}^{k+m+n},0)$.\\ 
We call $(\Phi,a)$ {\em a
reticular $t$-${\cal P}$-${\cal K}$-isomorphism from $F$ to $G$}. \\

We investigated the theory of function germs with respect to these equivalence 
relations in \cite{tPKfunct}. The main result is the following:
\begin{thm}\label{mthft:th}
Let $f$ be an unfolding of $f_0(x,y)\in {\mathfrak M}(r;k)$
and $F(x,y,t,u)\in {\mathfrak M}(r;k+m+n)$ be an unfolding of 
$f(x,y,u)\in {\mathfrak M}(r;k+n)$.  
Then the
following are all equivalent. \\ 
{\rm (1)} There exists a non-negative number $l$ such that 
$f_0$ is reticular 
${\cal K}$-$l$-determined and  $F$ is reticular
$t$-${\cal P}$-${\cal K}$-$l'$-transversal for $l'\geq lm+l+m+1$.\\
{\rm (2)} $F$ is reticular $t$-${\cal P}$-${\cal K}$-stable.\\
{\rm (3)} $F$ is reticular $t$-${\cal P}$-${\cal K}$-versal.\\ {\rm (4)} $F$ is reticular
$t$-${\cal P}$-${\cal K}$-infinitesimally versal.\\ 
{\rm (5)} $F$ is reticular
$t$-${\cal P}$-${\cal K}$-infinitesimally stable.\\ 
{\rm (6)} $F$ is reticular
$t$-${\cal P}$-${\cal K}$-homotopically stable.
\end{thm}
The classification list of reticular $t$-${\cal P}$-${\cal K}$-stable unfoldings 
in ${\mathfrak M}(r;k+1+n)$ with  $r=0,n\leq 5$ and $r=1,n\leq 3$ are 
given in \cite[p.201]{tPKfunct}.
\subsection{Contact regular $r$-cubic configurations}\label{conta:sec}
\quad 
We review the results given in \cite{retLeg}.
Let $J^1({\mathbb R}^n,{\mathbb R})$ be the $1$-jet bundle of functions 
in $n$-variables which may be considered as ${\mathbb R}^{2n+1}$ with 
a natural coordinate system $(q,z,p)=(q_1,\ldots,q_n,z,p_1,\ldots,p_n)$,
where $q$ be a coordinate system of ${\mathbb R}^n$.
We equip the contact structure on $J^1({\mathbb R}^n,{\mathbb R})$
defined by the canonical $1$-form $\theta=dz-\sum_{i=1}^np_idq_i$.
We have a natural projection 
$\pi:J^1({\mathbb R}^n,{\mathbb R})\rightarrow {\mathbb R}^n\times {\mathbb R}$
by $\pi(q,z,p)=(q,z)$.
\begin{defn}\label{config}{\rm
Let $w\in J^1({\mathbb R}^n,{\mathbb R})$ and $\{L_\sigma \}_{\sigma \subset I_r}$ be a family of $2^r$
 Legendrian submanifold germs on $(J^1({\mathbb R}^n,{\mathbb R}),w)$. 
Then $\{L_\sigma \}_{\sigma \subset I_r}$ is called 
{\em a contact regular $r$-cubic configuration} 
on $J^1({\mathbb R}^n,{\mathbb R})$
if there exists a 
contact embedding germ $C:(J^1({\mathbb R}^n,{\mathbb R}),0)
\rightarrow (J^1({\mathbb R}^n,{\mathbb R}),w)$
 such that 
$L_\sigma=C(L^0_\sigma)$ for all $\sigma \subset I_r$.
}\end{defn}
\begin{thm}{\rm (cf., \cite[Theorem 5.6]{retLeg})}
{\rm (1)} For any contact
 regular $r$-cubic configuration 
$\{L_\sigma \}_{\sigma \subset I_r}$ on 
$(J^1({\mathbb R}^n,{\mathbb R}),0)$, 
there exists a function germ $F\in {\mathfrak M}(r;k+n+1)$ which is a generating family of $\{L_\sigma \}_{\sigma \subset I_r}$. \\
{\rm (2)} For any C-non-degenerate function germ $F\in {\mathfrak M}(r;k+n+1)$, there exists a  contact
 regular $r$-cubic configuration  of which $F$ is a generating family. \\
{\rm (3)} Two  contact
 regular $r$-cubic configurations  are Legendrian equivalent if and only if their generating families are stably reticular ${\cal P}$-${\cal K}$-equivalent. 
\end{thm}
\section{Reticular Legendrian unfoldings}\label{RLeg:unfo}
\quad 
We consider the big $1$-jet bundle 
$J^1({\mathbb R}\times{\mathbb R}^n,{\mathbb R})$
and the canonical $1$-form $\Theta$ on that space.
Let $(t,q)=(t,q_1,\ldots,q_n)$ be the canonical coordinate system on 
${\mathbb R}\times{\mathbb R}^n$ and $(t,q,z,s,p)=(t,q_1,\ldots,q_n,z,s,p_1,\ldots,p_n)$ 
be the corresponding coordinate system on 
$J^1({\mathbb R}\times{\mathbb R}^n,{\mathbb R})$.
Then the canonical $1$-form $\Theta$ is given by 
\[ \Theta=dz-\sum_{i=1}^np_idq_i-sdt=\theta-sdt. \]

We recall that our purpose is the investigation of
one-parameter families of 
contact regular $r$-cubic configurations on 
$J^1({\mathbb R}^n,{\mathbb R})$ which defined by 
one-parameter families of 
contact embedding gems $(J^1({\mathbb R}^n,{\mathbb R}),0)
\rightarrow J^1({\mathbb R}^n,{\mathbb R})$ depending 
smoothly on $t\in ({\mathbb R},0)$.

Let $\{ L_{\sigma,t} \}_{\sigma\subset I_r,t\in ({\mathbb R},0)}$ be a family of contact regular $r$-cubic configurations on 
$J^1({\mathbb R}^n,{\mathbb R})$ defined by a family of contact embedding germs $C_t:
(J^1({\mathbb R}^n,{\mathbb R}),0)
\rightarrow J^1({\mathbb R}^n,{\mathbb R})$ depending 
smoothly on $t\in ({\mathbb R},0)$ such that $C_0(0)=0$ and $L_{\sigma,t} =
C_t(L^0_\sigma) $
for all $\sigma\subset I_r$ and $t\in ({\mathbb R},0)$.\\
Then we consider the following contact diffeomorphism germ $C$ on $(J^1({\mathbb R}\times {\mathbb R}^n,
{\mathbb R}),0)$:
\begin{lem}\label{C:lem}
For any family of contact embedding germs $C_t:
(J^1({\mathbb R}^n,{\mathbb R}),0)
\rightarrow J^1({\mathbb R}^n,{\mathbb R})\ (C_0(0)=0)$ depending 
smoothly on $t\in ({\mathbb R},0)$,
there exists a unique function germ $h$ on 
$(J^1({\mathbb R}\times {\mathbb R}^n,{\mathbb R}),0)$ such that 
the map germ $C:(J^1({\mathbb R}\times {\mathbb R}^n,{\mathbb R}),0)
\rightarrow (J^1({\mathbb R}\times {\mathbb R}^n,{\mathbb R}),0)$ defined by 
\[ C(t,q,z,s,p)=(t,q\circ C_t(q,z,p),z\circ C_t(q,z,p),
h(t,q,z,s,p),p\circ C_t(q,z,p))\]
is a contact diffeomorphism.
\end{lem}
{\em Proof}. 
We denote $C_t(q,z,p)=(q_t(q,z,p),z_t(z,q,p),p_t(q,z,p))$.
Since $C_t$ is a contact embedding germ for all $t\in ({\mathbb R},0)$,
there exists a function germ $\alpha(t,q,z,p)$ around zero with 
$\alpha(0)\neq 0$ such that 
$dz_t(q,z,p)-p_t(q,z,p)dq_t(q,$ $z,p)=\alpha(t,q,z,p)(dz-pdq)$
for all fixed t.
By the direct calculation of this equation, we have that
\[ 
 \frac{\partial z_t}{\partial z}-p_t
\frac{\partial q_t}{\partial z} = \alpha,
 \frac{\partial z_t}{\partial q}-p_t
\frac{\partial q_t}{\partial q}= -p \alpha,
 \frac{\partial z_t}{\partial p}-p_t
\frac{\partial q_t}{\partial p}=  0.\]
We also calculate $C^*(dz-pdq-sdt)$ by considering the above 
relations. Then we have that
\begin{eqnarray*}
 & & C^*(dz-pdq-sdt) \\
& = & dz_t(z,q,p)-p_t(q,z,p)dq_t(q,z,p)-
h(t,q,z,s,p)dt\\
& = & \alpha(t,z,q,p)dz-\alpha(t,z,q,p)pdq-
(\frac{\partial z_t}{\partial t}(q,z,p)-p_t(q,z,p)
\frac{\partial q_t}{\partial t}(q,z,p)-h(t,q,z,s,p))dt.
\end{eqnarray*}
In order to make $C$ a contact embedding, the function $h(t,q,z,s,p)$ 
is uniquely determined that:
\begin{equation}
h(t,q,z,s,p)=\frac{\partial z_t}{\partial t}(q,z,p)-p_t(q,z,p)
\frac{\partial q_t}{\partial t}(q,z,p)+\alpha(t,q,z,p)s. \label{hs:eqn}
\end{equation}
\hfill $\blacksquare$

\begin{defn}\label{C:dfn}{
Let $C$ be a contact diffeomorphism germ on
$(J^1({\mathbb R}\times {\mathbb R}^n,{\mathbb R}),0)$.
We say that $C$ is {\em a ${\cal P}$-contact diffeomorphism}
if $C$ has the form:
\begin{equation}
C(t,q,z,s,p)=(t,q_C(t,q,z,p),z_C(t,q,z,p),
h_C(t,q,z,s,p),p_C(t,q,z,p)).\label{tconta:eqn}
\end{equation}
}\end{defn}
We remark that a ${\cal P}$-contact diffeomorphism and the corresponding 
one-parameter family of contact embedding germs are uniquely defined by
each other.\\

We define that  
$\tilde{L}^0_\sigma=\{(t,q,z,s,p) \in J^1({\mathbb R}\times 
{\mathbb R}^n,{\mathbb R})|
q_\sigma=p_{I_r-\sigma}=q_{r+1}=\cdots=q_n=s=z=0,q_{I_r-\sigma}\geq 0 \}$
for $\sigma\subset I_r$
and ${\mathbb L}=\{(t,q,z,s,p) \in J^1({\mathbb R}\times 
{\mathbb R}^n,{\mathbb R})|
q_1p_1=\cdots=q_rp_r=q_{r+1}=\cdots=q_n=s=z=0,q_{I_r}\geq 0 \}$ be a 
representative as a germ 
of the union of $\tilde{L}^0_\sigma$ for all $\sigma\subset I_r$.

\begin{defn}{\rm 
We say that a map germ
${\cal L}:({\mathbb L},0)\rightarrow 
(J^1({\mathbb R}\times {\mathbb R}^n,{\mathbb R}),0)$ is
{\em a reticular Legendrian unfolding} if 
${\cal L}$ is the restriction of a ${\cal P}$-contact diffeomorphism.
We call $\{ {\cal L}(\tilde{L}^0_\sigma) \}_{\sigma\subset I_r}$ 
{\em the unfolded contact regular $r$-cubic configuration of } ${\cal L}$.
}\end{defn}

We note that:
Let $\{ \tilde{L}_{\sigma} \}_{\sigma\subset I_r}$ be
an unfolded contact regular $r$-cubic configuration
 associated with a one-parameter family of contact regular $r$-cubic configurations $\{ L_{\sigma,t} \}_{\sigma\subset I_r,t\in ({\mathbb R},0)}$.
Then there is the following relation between the wavefront 
$W_\sigma=\Pi(\tilde{L}_{\sigma})$ and 
the family of wavefronts $W_{\sigma,t}=\pi(L_{\sigma,t})$:
\[ W_\sigma=\bigcup_{t\in ({\mathbb R},0)}
\{t\}\times W_{\sigma,t} \ \ \ \ \mbox{ for all } \sigma\subset I_r.\]

In order to study bifurcations of wavefronts of 
 unfolded contact regular $r$-cubic configurations
we introduce the following equivalence relation.
Let $K,\Psi$ be contact diffeomorphism germs on 
$(J^1({\mathbb R}\times {\mathbb R}^n,{\mathbb R}),0)$.
We say that $K$ is {\em a ${\cal P}$-Legendrian equivalence} if $K$ has 
the form:
\begin{equation}
K(t,q,z,s,p)= (\phi_1(t),\phi_2(t,q,z),\phi_3(t,q,z),\phi_4(t,q,z,s,p),
\phi_5(t,q,z,s,p))\label{PLequi}.
\end{equation}
We say that $\Psi$ is {\em a reticular ${\cal P}$-diffeomorphism} if $\pi_t\circ \Psi$ depends only on $t$ and $\Psi$ preserves 
$\tilde{L}^0_\sigma$ for all $\sigma\subset I_r$. 

Let $\{ \tilde{L}^i_{\sigma} \}_{\sigma\subset I_r}(i=1,2)$
be unfolded contact regular $r$-cubic configurations on 
$(J^1({\mathbb R}\times {\mathbb R}^n,{\mathbb R}),0)$.
We say that they are 
{\em ${\cal P}$-Legendrian equivalent} 
if there exist a ${\cal P}$-contact diffeomorphism germ $K$
 such that 
$\tilde{L}^2_{\sigma}=K(\tilde{L}^1_{\sigma})$
for all $\sigma\subset I_r$.\vspace{2mm}

In order to understand the meaning of ${\cal P}$-Legendrian equivalence,
we observe the following: 
Let  $\{ \tilde{L}^i_{\sigma} \}_{\sigma\subset I_r}(i=1,2)$  be 
unfolded contact regular $r$-cubic configurations on 
$(J^1({\mathbb R}\times {\mathbb R}^n,{\mathbb R}),0)$
and 
$\{ L^i_{\sigma,t} \}_{\sigma\subset I_r,t\in ({\mathbb R},0)}$
be the corresponding one-parameter families of 
contact regular $r$-cubic configurations on 
$J^1({\mathbb R}^n,{\mathbb R})$
 respectively.
We take the smooth path germs $w_i:({\mathbb R},0)\rightarrow 
(J^1({\mathbb R}^n,{\mathbb R}),0)$ such that 
$\{ L^i_{\sigma,t} \}_{\sigma\subset I_r}$ are defined at $w_i(t)$ 
for $i=1,2$.
Suppose that there exists a ${\cal P}$-Legendrian equivalence $K$ from 
$\{ \tilde{L}^1_{\sigma} \}_{\sigma\subset I_r}$
to $\{ \tilde{L}^2_{\sigma} \}_{\sigma\subset I_r}$ of the form (\ref{PLequi}).
We set  $W_{\sigma,t}^i$ be the wavefront of $L^i_{\sigma,t}$ 
for $\sigma\subset I_r,\ t\in ({\mathbb R},0)$ and $i=1,2$.
We define the family of diffeomorphism $g_t:({\mathbb R}^n\times {\mathbb R},
\pi(w_1(t)))\rightarrow ({\mathbb R}^n\times {\mathbb R},
\pi(w_2(t)))$ by 
$g_t(q,z)=(\phi_2(t,q,z),\phi_3(t,q,z))$.
Then we have that $g_t(W^1_{\sigma,t})=W^1_{\sigma,\phi_1(t)}$ for all $\sigma\subset I_r,\ t\in ({\mathbb R},0)$.
\vspace{2mm}

We also define the equivalence relation among reticular Legendrian unfoldings.
Let ${\cal L}_i:({\mathbb L},0)\rightarrow 
(J^1({\mathbb R}\times {\mathbb R}^n,{\mathbb R}),0), (i=1,2)$ 
be reticular Legendrian unfoldings.
We say that ${\cal L}_1$ and ${\cal L}_2$ are 
{\em ${\cal P}$-Legendrian equivalent} 
if there exist a ${\cal P}$-Legendrian equivalence $K$ 
and a reticular ${\cal P}$-diffeomorphism 
$\Psi$ 
such that 
$K\circ {\cal L}_1={\cal L}_2\circ \Psi$.

We remark that two reticular Legendrian unfoldings are ${\cal P}$-Legendrian equivalent
if and only if the corresponding unfolded contact regular $r$-cubic configurations are ${\cal P}$-Legendrian equivalent.
\vspace{2mm}

By the same proof of Lemma 5.3 in \cite{retLeg}, we have the following:
\begin{lem}\label{exleg:lm}
Let $\{\tilde{L}_\sigma \}_{\sigma \subset I_r}$ be an 
unfolded contact regular $r$-cubic configuration on 
$(J^1({\mathbb R}\times {\mathbb R}^n,{\mathbb R}),0)$.
Then there exists 
a ${\cal P}$-contact diffeomorphism germ $C$ on $(J^1({\mathbb R}\times {\mathbb R}^n,{\mathbb R}),0)$
such that $C$ defines
$\{\tilde{L}_\sigma \}_{\sigma \subset I_r}$  and preserves the canonical 
$1$-form.
\end{lem}
By this lemma we may 
assume that all reticular Legendrian unfoldings (and all unfolded contact regular $r$-cubic configurations) 
are defined by ${\cal P}$-contact diffeomorphism germs which 
preserve the canonical $1$-form.
 \vspace{2mm}

We can construct generating families of reticular Legendrian unfoldings.
A function germ $F(x,y,t,q,z)\in{\mathfrak M}(r;k+1+n+1)$ is said to be 
{\em ${\cal P}$-$C$-non-degenerate} if 
$\frac{\partial F}{\partial x}(0)=\frac{\partial F}{\partial y}(0)=0$ and 
$x,t,F,\frac{\partial F}{\partial x},
\frac{\partial F}{\partial y}$
are independent on $({\mathbb H}^k\times {\mathbb R}^{k+1+n+1},0)$.\\

A ${\cal P}$-$C$-non-degenerate function germ 
$F(x,y,t,q,z)\in {\mathfrak M}(r;k+1+n+1)$
is called {\em a generating family} of 
a reticular Legendrian unfoldings ${\cal L}$ if
\begin{eqnarray*}
{\cal L}(\tilde{L}^0_{\sigma}) =
\{ (t,q,z,\frac{\partial F}{\partial t}/(-\frac{\partial F}{\partial z}),
\frac{\partial F}{\partial q}/(-\frac{\partial F}{\partial z}))\in
(J^1({\mathbb R}\times {\mathbb R}^n,{\mathbb R}),0)|\hspace{1cm}\\
\hspace{2cm}
x_\sigma=F=\frac{\partial F}{\partial x_{I_r-\sigma}}=
\frac{\partial F}{\partial y}=0,x_{I_r-\sigma}\geq 0\} 
\mbox{ for all }\sigma\subset I_r.
\end{eqnarray*}
We remark that for a ${\cal P}$-$C$-non-degenerate function germ 
$F(x,y,t,q,z)$, 
the function germ $F(\cdot,\cdot,t,\cdot,\cdot)$ is $C$-non-degenerate (see \cite[p.111]{retLeg}). 
\vspace{2mm}

\begin{lem}\label{g.function of rLu:lem}
Let $C$ be a ${\cal P}$-contact diffeomorphism germ on 
$(J^1({\mathbb R}\times {\mathbb R}^n,{\mathbb R}),0)$ 
which preserves the canonical $1$-form.
If 
the map germ
\[ (T,Q,Z,S,P)\rightarrow (T,Q,Z,s_C(T,Q,Z,S,P),p_C(T,Q,Z,S,P))\]
is diffeomorphism at $0$,
then there exists a function germ $H(T,Q,p)\in 
{\mathfrak M}(1+n+n)^2$ such that 
the canonical relation $P_C$ associated with $C$ has the form:
\begin{equation}
P_C=\{ (T,Q,Z,-\frac{\partial H}{\partial T}(T,Q,p)+s,
-\frac{\partial H}{\partial Q},\hspace{3cm}\\
T,-\frac{\partial H}{\partial p},
H-\langle \frac{\partial H}{\partial p},p\rangle +Z,
s,p)\},\label{cano:eqn}
\end{equation}
and the function germ $F\in {\mathfrak M}(r;n+1+n+1)$ defined by
$F(x,y,t,q,z)=-z+H(t,x,0,y)+\langle y,q\rangle$
is a generating family of the reticular Legendrian unfolding $C|_{{\mathbb L}}$.
\end{lem}
{\em Proof}.
We have that $dz-sdt-pdq=dZ-SdT-PdQ$ on $P_C$.
It follows that $d(z-Z)=sdt+pdq-SdT-PdQ$ and 
$d(z-Z+st+pq)=-tds-qdp-SdT-PdQ$.
Then there exists a function germ $H'(T,Q,s,p)\in 
{\mathfrak M}(1+n+1+1+n)^2$ such that 
\[
z-Z-st-pq=H'(T,Q,s,p),\ t=-\frac{\partial H'}{\partial s},\ 
q=-\frac{\partial H'}{\partial p},
S=-\frac{\partial H'}{\partial T},\
P=-\frac{\partial H'}{\partial Q}
\mbox{ on } P_C. \]
Since $t=T=-\frac{\partial H'}{\partial s}$ on $P_C$, 
we have that $H'(T,Q,s,p)=H(T,Q,p)-Ts$
for some $H(T,Q,p)\in {\mathfrak M}(1+n+n)^2$.
Then we have that
\[ z-Z-Ts- \langle -\frac{\partial H}{\partial p},p\rangle=H(T,Q,p)-Ts.\]
It follows that
\[ z=H(T,Q,p)-\langle \frac{\partial H}{\partial p},p\rangle+Z. \]
Then we have the required form of $P_C$.
By the direct calculation with the form $P_C$, 
we have that $F$ is a generating family
of $C|_{{\mathbb L}}$.\hfill $\blacksquare$\\

We have the following theorem which gives the relations 
between reticular Legendrian unfoldings and their generating families.
\begin{thm}\label{UCgf:th}
{\rm (1)} For any reticular Legendrian unfolding ${\cal L}:({\mathbb L},0)
\rightarrow (J^1({\mathbb R}\times {\mathbb R}^n,{\mathbb R}),0)$, 
there exists a function germ $F(x,y,t,q,z)\in{\mathfrak M}(r;k+1+n+1)$ 
which is a generating family of ${\cal L}$.\\
{\rm (2)} For any ${\cal P}$-$C$-non-degenerate function germ 
$F(x,y,t,q,z)\in{\mathfrak M}(r;k+1+n+1)$ with
$\frac{\partial F}{\partial t}(0)=\frac{\partial F}{\partial q}(0)=0$, there exists a reticular
 Legendrian unfolding ${\cal L}:({\mathbb L},0)\rightarrow 
(J^1({\mathbb R}\times {\mathbb R}^n,{\mathbb R}),0)$ 
of which $F$
is a generating family.\\
{\rm (3)} Two reticular Legendrian unfolding are ${\cal P}$-Legendrian 
equivalent if and only if their generating families are stably reticular 
$t$-${\cal P}$-${\cal K}$-equivalent.
\end{thm}
This theorem is proved by analogous methods of \cite{retLag}, \cite{retLeg}.
We give the sketch of the proof.
(1) Let $C$ be a ${\cal P}$-contact diffeomorphism germ on 
$(J^1({\mathbb R}\times {\mathbb R}^n,{\mathbb R}),0)$ 
such that $C|_{{\mathbb L}}={\cal L}$.
We may assume that $C^*\Theta=\Theta$.
By taking a ${\cal P}$-Legendrian equivalence of ${\cal L}$,
we may assume that the canonical relation $P_C$ associated with $C$ has the
form (\ref{cano:eqn}) for the function germ $H\in {\mathfrak M}(1+n+n)^2$.
Then the function germ $F(x,y,t,q,z)\in {\mathfrak M}(r;n+1+n+1)$ defined
 by
\[F(x,y,t,q,z)=-z+H(t,x_1,\ldots,x_r,0,y)+\langle y,q \rangle \]
is a generating family of ${\cal L}$.\\
(2) Let a ${\cal P}$-$C$-non-degenerate function germ $F(x,y,t,q,z)\in{\mathfrak M}(r;k+1+n+1)$ with
$\frac{\partial F}{\partial t}(0)=\frac{\partial F}{\partial q}(0)=0$ be given.
By \cite[Lemma 2.1]{retLeg},
we may assume that $F$ has the form $F(x,y,t,q,z)=-z+F_0(x,y,t,q)$
for some $F_0\in {\mathfrak M}(r;k+1+n)$.
Choose an $(n-r)\times k$-matrix $A$ and an $(n-r)\times n$-matrix $B$ 
such that the matrix
\begin{equation}
 \left( 
\begin{array}{ccc}
\displaystyle{\frac{\partial^2 F_0}{\partial x\partial y} } & 
\displaystyle{\frac{\partial^2 F_0}{\partial x\partial q} } &
\displaystyle{\frac{\partial^2 F_0}{\partial x\partial t} } \\
\displaystyle{\frac{\partial^2 F_0}{\partial y\partial y} } & 
\displaystyle{\frac{\partial^2 F_0}{\partial y\partial q} } & 
\displaystyle{\frac{\partial^2 F_0}{\partial y\partial t} } \\
A & B & 0 \\
0 & 0 & 1 
\end{array}
\right)_0\mbox{ is invertible. } 
\end{equation}
Let $F'\in {\mathfrak M}(r+k+1+n+1)$ be a function germ
which is obtained by an extension the source space of $F$ to 
$({\mathbb R}^{r+k+1+n+1},0)$.
Define the function $G(S,Q,y,t,q,z)\in  
{\mathfrak M}(n+1+1+k+1+n+1)$ by that 
\[ G(Q,Z,S,y,t,q,z)=-z+F'(Q_1,\ldots,Q_r,y,t,q)+\hspace{4cm}\] 
\[\hspace{1cm} (Q_{r+1},\ldots,Q_n) A 
\left(
\begin{array}{c}
y_1\\ \vdots \\ y_k 
\end{array}\right)
+(Q_{r+1},\ldots,Q_n) B
\left(
\begin{array}{c}
q_1\\ \vdots \\ q_n 
\end{array}
\right)+St.
\]
Then $G$ is a generating family of the canonical relation $P_C$ associated
with some ${\cal P}$-contact diffeomorphism germ $C$.
The function germ $F$ is a generating family of the reticular Legendrian
unfolding $C|_{{\mathbb L}}$.\\
(3) We need only to prove that: {\em If $F_1,F_2\in {\mathfrak M}(r;k+1+n)$ are 
generating families of the same reticular Legendrian unfolding, then they are
reticular
$t$-${\cal P}$-${\cal K}$-equivalent.}\\
We may reduce that $F_i$ has the form $F_i(x,y,t,q,z)=-z+F_i^0(x,y,t,q)$ for
$F^0_i\in {\mathfrak M}(r;k+1+n),\ i=1,2$.
Then $F^0_1$ and $F^0_2$ are generating families of the same reticular Lagrangian map in the sense of \cite{retLag}.
By \cite[p.587 the assertion (3)]{retLag},
there exists a reticular ${\cal R}$-equivalence from $F^0_2$ to $F^0_1$ of the 
form:
\[ F^0_1(x,y,t,q)=F^0_2(x\phi_1(x,y,t,q),\phi_2(x,y,t,q),t,q).\]
This means that $F_1$ and $F_2$ are reticular
$t$-${\cal P}$-${\cal K}$-equivalent.\hfill $\blacksquare$

\section{Stabilities}\label{map:sec}
\quad  
In this section we shall define several stabilities of reticular 
Legendrian unfoldings and prove that they and the 
stabilities of corresponding generating families are all equivalent.

Let $U$ be an open set in $J^1({\mathbb R}\times {\mathbb R}^n,{\mathbb R})$.
We consider contact diffeomorphism germs on 
$(J^1({\mathbb R}\times {\mathbb R}^n,{\mathbb R}),0)$ 
and contact embeddings
from $U$ 
to $J^1({\mathbb R}\times {\mathbb R}^n,{\mathbb R})$. 
Let $(T,Q,S,Z,P)$ and $(t,q,z,s,p)$ be canonical coordinates of the source
space and the target space respectively. 
We define the following notations:\\
$\imath:(J^1({\mathbb R}\times {\mathbb R}^n,{\mathbb R})\cap \{Z=0 \},0)\rightarrow 
(J^1({\mathbb R}\times {\mathbb R}^n,{\mathbb R}),0)$ be the inclusion map on the source space,
\begin{eqnarray*}
C_T(J^1({\mathbb R}\times {\mathbb R}^n,{\mathbb R}),0) & = &
\{ C| C \mbox{ is a ${\cal P}$-contact diffeomorphism  germ on }(J^1({\mathbb R}\times {\mathbb R}^n,{\mathbb R}),0) \}, \\
C_T^\Theta (J^1({\mathbb R}\times {\mathbb R}^n,{\mathbb R}),0)
& = & \{ C\in C_T(J^1({{\mathbb R}\times \mathbb R}^n
,{\mathbb R}),0)|\ C^*\Theta=\Theta \}, \\
C_T^{Z} (J^1({\mathbb R}\times {\mathbb R}^n,{\mathbb R}),0)
 & = & \{ C\circ\imath\ |C \in 
C_T(J^1({\mathbb R}\times {\mathbb R}^n,{\mathbb R}),0) \},\\
C_T^{\Theta,Z} (J^1({\mathbb R}\times {\mathbb R}^n,{\mathbb R}),0)
 & = & \{ C\circ\imath\ |
C \in C_T^\Theta (J^1({\mathbb R}\times {\mathbb R}^n,{\mathbb R}),0) \}.
 \end{eqnarray*}
 
Let $V=U\cap \{Z=0\}$ and $\tilde{\imath}:V\rightarrow U$ be the 
inclusion map.
\begin{eqnarray*}
C_T(U,J^1({\mathbb R}\times {\mathbb R}^n,{\mathbb R})) & = & \{ 
\tilde{C}:U\rightarrow 
J^1({\mathbb R}\times {\mathbb R}^n,{\mathbb R})| \\
 & & \hspace{1cm} \tilde{C} \mbox{ is a contact embedding of the form (\ref{tconta:eqn})}\},\\
 C_T^\Theta (U,J^1({\mathbb R}\times {\mathbb R}^n,{\mathbb R}))  & = & \{ \tilde{C}\in 
C_T(U,J^1({\mathbb R}\times {\mathbb R}^n,{\mathbb R}))\ |
\tilde{C}^*\Theta=\Theta \},\\
 C_T^Z (V,J^1({\mathbb R}\times {\mathbb R}^n,{\mathbb R})) & = & \{ \tilde{C}\circ
\tilde{\imath}\ |\tilde{C}\in C_T(U,J^1({\mathbb R}\times {\mathbb R}^n,{\mathbb R}) )\},\\
 C_T^{\Theta,Z} (V,J^1({\mathbb R}\times {\mathbb R}^n,{\mathbb R})) & = & 
 \{ \tilde{C}\circ
\tilde{\imath}\ |\tilde{C}\in C_T^\Theta (U,J^1({\mathbb R}\times {\mathbb R}^n,{\mathbb R})) \}.
\end{eqnarray*}

\begin{defn}{\rm
{\bf Stability}: 
We say that a reticular Legendrian unfolding ${\cal L}$ is {\em stable} if the 
following condition holds:
Let $C_0\in C_T(J^1({\mathbb R}\times {\mathbb R}^n,{\mathbb R}),0))$ 
be a ${\cal P}$-contact diffeomorphism germ 
such that $C_0|_{{\mathbb L}}={\cal L}$ and 
$\tilde{C}_0\in C_T(U,J^1({\mathbb R}\times {\mathbb R}^n,{\mathbb R}))$ be
a representative of $C_0$.
Then there exists an open neighborhood $N_{\tilde{C}_0}$ of $\tilde{C}_0$ 
in $C^\infty$-topology
such that for any $\tilde{C}\in N_{\tilde{C}_0}$,
there exists a point $w_0=(T^0,0,\ldots,0,P^0_{r+1},\ldots,P^0_n)\in U$ such that 
the reticular Legendrian unfolding ${\cal L}'_{w_0}$ and 
${\cal L}$ are 
${\cal P}$-Legendrian equivalent,
where the reticular Legendrian unfolding ${\cal L}'_{w_0}$ is defined by
\[ x=(T,Q,Z,S,P)\mapsto \tilde{C}(w_0+x)-\tilde{C}(w_0)+(0,0,P^0_{r+1}Q_{r+1}+\cdots +P^0_nQ_n,0,0). \] 
{\bf Homotopical stability}:
A one-parameter family
of ${\cal P}$-contact diffeomorphism germs  $\bar{C}:(J^1({\mathbb R}\times {\mathbb R}^n,{\mathbb R})
\times {\mathbb R},(0,0))\rightarrow
(J^1({\mathbb R}\times {\mathbb R}^n,{\mathbb R}),0)((T,Q,Z,S,P,\tau)
\mapsto C_\tau(T,Q,Z,S,P))$
is called a {\em reticular ${\cal P}$-contact 
  deformation} of ${\cal L}$ if $C_0|_{{\mathbb L}}={\cal L}$.
A map germ $\bar{\Psi}:(J^1({\mathbb R}\times {\mathbb R}^n,{\mathbb R}) \times 
{\mathbb R},(0,0)) \rightarrow 
(J^1({\mathbb R}\times {\mathbb R}^n,{\mathbb R}),0)((T,Q,Z,S,P,\tau)
\mapsto \Psi_\tau(T,Q,Z,S,P))$
is called a {\em one-parameter
  deformation of reticular ${\cal P}$-diffeomorphisms} if 
$\Psi_0=id_{J^1({\mathbb R}\times {\mathbb R}^n,{\mathbb R})}$ and 
$\Psi_t$ 
 is a reticular ${\cal P}$-diffeomorphism for all $t$ around $0$. 
We say that a reticular Legendrian unfolding ${\cal L}$ is {\em
  homotopically stable} if for any reticular ${\cal P}$-Legendrian deformation 
$\bar{C}=\{ C_\tau\}$ of ${\cal L}$,
there exist 
 one-parameter family of ${\cal P}$-Legendrian 
  equivalences $\bar{K}=\{ K_\tau \}$ with $K_0=id_{J^1({\mathbb R}\times {\mathbb R}^n,{\mathbb R})}$ and 
a one-parameter deformation of
reticular ${\cal P}$-diffeomorphisms $\bar{\Psi}=\{ \Psi_\tau \}$
such that $C_\tau=K_\tau\circ C_0\circ \Psi_\tau$ for $t$ around
$0$.\\
{\bf Infinitesimal stability}:
Let $C\in C_T(J^1({\mathbb R}\times {\mathbb R}^n,{\mathbb R}),0)$
be a ${\cal P}$-contact diffeomorphism germ.
We say that a vector field $v$ on $(J^1({\mathbb R}\times {\mathbb R}^n,{\mathbb R}),0)$ along $C$ 
is {\em an infinitesimal ${\cal P}$-contact transformation} of $C$ if 
there exists a  reticular ${\cal P}$-Legendrian deformation 
$\bar{C}=\{C_\tau\}$ on $(J^1({\mathbb R}\times {\mathbb R}^n,{\mathbb R}),0)$ 
such that $C_0=C$  and
$\frac{dC_\tau}{d\tau}|_{\tau =0}=v$.
We say that a vector field $\xi$ on $(J^1({\mathbb R}\times {\mathbb R}^n,{\mathbb R}),0)$ is 
{\em infinitesimal reticular ${\cal P}$-diffeomorphism} if there exists 
a one-parameter deformation of
reticular ${\cal P}$-diffeomorphisms  $\bar{\Psi}=\{ \Psi_\tau \}$ 
such that  $\frac{d\Psi_\tau}{d\tau }|_{\tau =0}=\xi$.
We say that a vector field $\eta$ on  $(J^1({\mathbb R}\times {\mathbb R}^n,{\mathbb R}),0)$ is
{\em infinitesimal ${\cal P}$-Legendrian equivalence} if 
there exists a
 one-parameter family of ${\cal P}$-Legendrian equivalences $\bar{K}=\{K_\tau\}$ such that 
$K_0=id_{J^1({\mathbb R}\times {\mathbb R}^n,{\mathbb R})}$
and $\frac{dK_\tau}{d\tau}|_{\tau =0}=\eta$.
We say that a reticular
Legendrian unfolding ${\cal L}$ is  
 {\em infinitesimally stable} if for any
extension $C\in C_T(J^1({\mathbb R}\times {\mathbb R}^n,{\mathbb R}),0)$ 
of ${\cal L}$ and any infinitesimal ${\cal P}$-contact transformation $v$ of $C$,
there exists an infinitesimal reticular ${\cal P}$-diffeomorphism $\xi$ and 
an infinitesimal ${\cal P}$-Legendrian equivalence $\eta$ such that $v=C_*\xi+\eta\circ C$.
}
\end{defn}

We may take an extension of a reticular Legendrian unfolding ${\cal L}$ by an element of $C^\Theta_T (J^1({\mathbb R}\times {\mathbb R}^n,{\mathbb R}),0)$ by
Lemma \ref{exleg:lm}.
Then as the remark after the definition of the stability of reticular 
Legendrian maps in \cite[p.121]{retLeg}, 
we may consider the following other definitions of stabilities of reticular Legendrian unfoldings: 
(1) The definition given by replacing 
$C_T(J^1({\mathbb R}\times {\mathbb R}^n,{\mathbb R}),0))$ and
$C_T(U,J^1({\mathbb R}\times {\mathbb R}^n,{\mathbb R}))$ to
$C_T^\Theta(J^1({\mathbb R}\times {\mathbb R}^n,{\mathbb R}),0))$ and 
$C_T^\Theta(U,J^1({\mathbb R}\times {\mathbb R}^n,{\mathbb R}))$
of original definition respectively.
(2) The definition given by replacing to $C_T^Z(J^1({\mathbb R}\times {\mathbb R}^n,{\mathbb R}),0))$ and
$C_T^Z(V,J^1({\mathbb R}\times {\mathbb R}^n,{\mathbb R}))$ respectively. 
(3) The definition given by replacing to 
$C_T^{\Theta,Z}(J^1({\mathbb R}\times {\mathbb R}^n,{\mathbb R}),0))$ and
$C_T^{\Theta,Z}(V,J^1({\mathbb R}\times {\mathbb R}^n,{\mathbb R}))$ 
respectively, where $V=U\cap \{Z=0\}$.\vspace{2mm}

Then we have the following lemma which is proved by the same method of the 
proof of Lemma 7.2 in \cite{retLeg}
\begin{lem}\label{sta:lm}
The original definition and other three definitions of stabilities of reticular Legendrian unfoldings are all equivalent.
\end{lem}
By this lemma, we may choose an extension of a reticular Legendrian unfolding 
from among all of 
 $C_T(J^1({\mathbb R}\times {\mathbb R}^n,{\mathbb R}),0))$, 
$C_T^\Theta(J^1({\mathbb R}\times {\mathbb R}^n,{\mathbb R}),0))$,
 $C_T^Z(J^1({\mathbb R}\times {\mathbb R}^n,{\mathbb R}),0))$, and 
$C_T^{\Theta,Z}(J^1({\mathbb R}\times {\mathbb R}^n,{\mathbb R}),0))$.\\

We say that a function germ $H$ on $(J^1({\mathbb R}\times {\mathbb R}^n,{\mathbb R}),0)$ is {\em ${\cal P}$-fiber preserving} if 
$H$ has the form $H(t,q,z,s,p)=\sum_{j=1}^nh_j(t,q,z)p_j+h_0(t,q,z)+a(t)s$.
\vspace{2mm}
\begin{lem}\label{infsta:t-Leglem}
Let $C\in C_T(J^1({\mathbb R}\times {\mathbb R}^n,{\mathbb R}),0)$.
Then the following hold:
{\rm (1)} A vector field germ $v$ on $(J^1({\mathbb R}\times {\mathbb R}^n,{\mathbb R}),0)$ along $C$
 is an infinitesimal ${\cal P}$-contact transformation of $C$ if and only if 
there exists a function germ $f$ on 
$(J^1({\mathbb R}\times {\mathbb R}^n,{\mathbb R}),0)$ 
such that $f$ does not depend on $s$ and $v=X_f\circ C$.\\
{\rm (2)} A vector field germ $\eta$ on 
$(J^1({\mathbb R}\times {\mathbb R}^n,{\mathbb R}),0)$ is
an infinitesimal ${\cal P}$-Legendrian equivalence
 if and only if there exists a ${\cal P}$-fiber preserving function germ
$H$ on $(J^1({\mathbb R}\times {\mathbb R}^n,{\mathbb R}),0)$ 
such that $\eta=X_H$.\\
{\rm (3)} A vector field $\xi$ on 
$(J^1({\mathbb R}\times {\mathbb R}^n,{\mathbb R}),0)$ is
an infinitesimal reticular ${\cal P}$-diffeomorphism  if and only if 
there exists a function germ $g\in B$
such that $\xi=X_g$, where $B=\langle q_1p_1,\ldots,q_rp_r,
q_{r+1},\ldots,q_n,z\rangle_{{\cal E}_{t,q,z,p}}
+\langle s \rangle_{{\cal E}_t}$
\end{lem}
{\em Proof}. 
A vector field $X$ on 
$(J^1({\mathbb R}\times {\mathbb R}^n,{\mathbb R}),0)$
is a contact Hamiltonian vector field if and only 
there exists a function germ $f$ on $(J^1({\mathbb R}\times {\mathbb R}^n,{\mathbb R}),0)$ such that $X=X_f$, that is
\begin{eqnarray*}
X=\sum_{i=1}^n(\frac{\partial f}{\partial q_i}-p_i
\frac{\partial f}{\partial z})\frac{\partial }{\partial p_i}+
(\frac{\partial f}{\partial t}+s\frac{\partial f}{\partial z})
\frac{\partial }{\partial s}\hspace{3cm}\\
-
\sum_{i=1}^n\frac{\partial f}{\partial p_i}\frac{\partial }{\partial q_i}-
\frac{\partial f}{\partial s}\frac{\partial }{\partial t}+
(f-\sum_{i=1}^np_i\frac{\partial f}{\partial p_i}-s
\frac{\partial f}{\partial s})\frac{\partial }{\partial z},
\end{eqnarray*}
(for definition, see \cite[p.145]{ishikawa2}).
(1) A vector field  $v$ is an  infinitesimal ${\cal P}$-contact transformation of $C$
if and only if $v=H_f\circ C$  and 
$\frac{\partial f}{\partial s}=0$.
This holds if and only  if
 $v=H_f\circ C$ and $f$ does not depend on $s$.\\
(2)  A vector field germ $\eta$ is an infinitesimal ${\cal P}$-Legendrian equivalence
 if and only if there exists a fiber preserving function germ
$H$ such that $\eta=X_H$ by 
and 
$\frac{\partial H}{\partial s}=a(t)$ for some function germ $a(t)$.
This holds if and only if $\eta=X_H$ and $H$ is a ${\cal P}$-fiber preserving function.\\
(3) A vector field $\xi$ is
an infinitesimal reticular ${\cal P}$-diffeomorphism of $C$ if and only if 
there exists a function germ $g\in \langle q_1p_1,\ldots,q_rp_r,
q_{r+1},\ldots,q_n,z,s\rangle_{{\cal E}_{J^1({\mathbb R}\times {\mathbb R}^n,{\mathbb R})}}$
such that $\xi=X_g$ 
since $X_g$ is tangent to $\tilde{L}^0_\sigma$ for all 
$\sigma \subset I_r$,
and  $\frac{\partial g}{\partial s}=a(t)$,
this holds if and only if 
$\xi=X_g$ and $g\in B$.
\hfill $\blacksquare$\\

Let $U$ be a neighborhood of $0$ in
$J^1({\mathbb R}\times {\mathbb R}^n,{\mathbb R})$.
We define: 
\begin{eqnarray*}
J^{l}_{C^\Theta_T}(U,J^1({\mathbb R}\times {\mathbb R}^n,
{\mathbb R}))=\{
j^lC(w_0)\in J^l(U,J^1({\mathbb R}\times {\mathbb R}^n,
{\mathbb R}))|\hspace{4cm}\\
\hspace{0.5cm}C:(U,w_0)\rightarrow J^1({\mathbb R}\times {\mathbb R}^n,
{\mathbb R})\mbox{ is a ${\cal P}$-contact embedding germ which 
preserves }
\Theta \}.
\end{eqnarray*}
\begin{thm}[${\cal P}$-Contact transversality theorem]\label{t_contra:th}
Let $Q_i,i=1,2,\ldots$ are submanifolds of 
$J^{l}_{C^\Theta_T}(U,J^1({\mathbb R}\times {\mathbb R}^n,
{\mathbb R}))$.
Then the set 
\[ T=\{ C\in C^\Theta_T(U,
J^1({\mathbb R}\times {\mathbb R}^n,{\mathbb R}))|
j^lC \mbox{ is transversal to } Q_i \mbox{ for all } i\in {\mathbb N}\}
\]
is a residual set in 
$C^\Theta_T(U,
J^1({\mathbb R}\times {\mathbb R}^n,{\mathbb R}))$
\end{thm}
This is proved by an almost parallel method of \cite[Theorem 6.4]{generic}.
\vspace{2mm}

We denote the ring ${\cal E}(1+n+n)$ 
on the coordinates $(t,q,p)$
by ${\cal E}_{t,q,p}$ and denote other notations analogously.
\begin{thm}\label{staLeg:th}
Let ${\cal L}$ 
be a reticular Legendrian unfolding with a 
generating family $F(x,y,t,q,z)$.
Then the following are all equivalent.\\
{\rm (u)} $F$ is a reticular $t$-${\cal P}$-${\cal K}$-stable unfolding of 
$F|_{t=0}$.\\
{\rm (hs)} ${\cal L}$ is homotopically stable.\\
{\rm (is)} ${\cal L}$ is infinitesimally stable.\\ 
{\rm (a)} 
${\cal E}_{t,q,p}=
B_0+
\langle 1,p_1\circ C',\ldots,p_n\circ C'\rangle_{(\Pi\circ C')^*{\cal E}_{t,q,z}}+
\langle s\circ C'\rangle_{{\cal E}_t}$,
where $C'=C|_{z=s=0}$ and $B_0=\langle q_1p_1,\ldots,q_rp_r,
q_{r+1},\ldots,q_n\rangle_{{\cal E}_{t,q,p}}$.
\end{thm}
We remark that sufficiently near reticular Legendrian unfoldings of 
stable one are stable by the condition (a).\\
{\em Proof}. (u)$\Rightarrow$(hs):
Let a reticular ${\cal P}$-Legendrian deformation 
$\bar{C}=\{ C_\tau\}$ of ${\cal L}$ be given.
The homotopically stability of reticular Legendrian unfoldings is invariant under 
${\cal P}$-Legendrian equivalences,
we may assume that 
the map germs 
\[ (T,Q,Z,S,P)\rightarrow (T,Q,Z,s\circ C_\tau(T,Q,P),p\circ C_\tau(T,Q,P))\] 
are diffeomorphisms for all $\tau$.
By Lemma \ref{g.function of rLu:lem}, there exists a one-parameter family 
$H_\tau(T,Q,p)\in {\mathfrak M}(1+n+n)^2$ depending smoothly on 
$\tau \in ({\mathbb R},0)$ such that the canonical relations $P_{C_\tau}$
associated with $C_\tau$ has the form:
\[
P_{C_\tau}=\{ (T,Q,Z,-\frac{\partial H_\tau}{\partial T}(T,Q,p)+s,
-\frac{\partial H_\tau}{\partial Q},
T,
-\frac{\partial H_\tau}{\partial p},
H_\tau-\langle \frac{\partial H_\tau}{\partial p},p\rangle +Z,
s,p)\}.
\]
Then the function germs $F_\tau\in {\mathfrak M}(r;n+1+n+1)$ defined by
\[
F_\tau(x,y,t,q,z)=-z+H_\tau(t,x,0,y)+
\langle y,q\rangle \]
are generating families of reticular Legendrian unfoldings ${\cal L}_\tau:=C_\tau|_{{\mathbb L}}$ for $\tau\in ({\mathbb R},0)$.
Since $F_0$ is a reticular $t$-${\cal P}$-${\cal K}$-stable unfolding of 
$F|_{t=0}$,
it follows that $F_0$ is a reticular $t$-${\cal P}$-${\cal K}$-homotopically 
stable unfolding of 
$F|_{t=0}$  by Theorem \ref{mthft:th}.
Therefore there exists a one-parameter family of reticular 
$t$-${\cal P}$-${\cal K}$-isomorphism 
from $F_\tau$ to $F_0$ depending smoothly on $\tau$.
This means that there exists a one-parameter family of ${\cal P}$-Legendrian equivalences 
$K_\tau$  depending smoothly on $\tau$ such that 
\[ C_\tau(L^0_\sigma)=K_\tau \circ {\cal L}(L^0_\sigma)
\mbox{ for all } \sigma\subset I_r,\ \ \tau\in ({\mathbb R},0).\]
Then the map germs $\Psi_\tau:=C_0^{-1}\circ K_\tau^{-1}\circ C_\tau$ give
 a one-parameter deformation of 
reticular ${\cal P}$-diffeomorphisms on 
$(J^1({\mathbb R}\times {\mathbb R}^n,{\mathbb R}),0)$ and 
we have  that $C_\tau=K_\tau \circ C_0\circ \Psi_\tau$.
This means that ${\cal L}$ is homotopically stable.\\
(hs)$\Rightarrow$(is):
Let $C\in C_T(J^1({\mathbb R}\times {\mathbb R}^n,{\mathbb R}),0)$ be 
an extension of ${\cal L}$ and
$v$ be an infinitesimal ${\cal P}$-contact transformation of $C$.
Then there exists a  reticular ${\cal P}$-Legendrian deformation 
$\bar{C}=\{C_\tau\}$ of $C$ such that $v=\frac{d C_\tau}{d \tau}|_{\tau=0}$.
Then there exist a one-parameter of ${\cal P}$-Legendrian 
 equivalences  $\bar{K}=\{ K_\tau \}$ and 
a 
 one-parameter
  deformation of reticular ${\cal P}$-diffeomorphisms $\bar{\Psi}=\{\Psi_\tau\}$
such that $C_\tau=K_\tau\circ C_0\circ \Psi_\tau$ for $\tau\in
({\mathbb R},0)$.
Then we have that 
\[ v=\frac{d C_\tau}{d \tau}|_{\tau=0}=
\frac{d K_\tau}{d \tau}|_{\tau=0}\circ C_0+(C_0)_*
(\frac{d \Psi_\tau}{d \tau}|_{\tau=0}).\]
(is)$\Rightarrow $(a):
Let a function germ $f\in {\cal E}_{t,q,p}$ be given.
We define the function germ $f'$ on 
$(J^1({\mathbb R}\times {\mathbb R}^n,{\mathbb R}),0)$ by 
$f'(t,q,z,s,p)=f\circ \pi_{T,Q,P}\circ C^{-1}(t,q,z,a(t,q,z,p),p)$,
where the function germ $a$ is defined by $\pi_S\circ C^{-1}(t,q,z,a(t,q,z,p),p)\equiv 0$.
This equation can be solved by (\ref{hs:eqn}).
Since $f'$ does not depend on $s$, it follows that 
$X_{f'}\circ C$ is an infinitesimal ${\cal P}$-contact transformation of $C$.
Therefore there exist an infinitesimal ${\cal P}$-Legendrian equivalence $\eta$ and 
an infinitesimal reticular ${\cal P}$-diffeomorphism $\xi$ such that 
$X_{f'}\circ C=C_*\xi +\eta\circ C$.
By Lemma \ref{infsta:t-Leglem}, there exist a ${\cal P}$-fiber preserving function germ $H$
on $(J^1({\mathbb R}\times {\mathbb R}^n,{\mathbb R}),0)$
and $g\in B$ such that $\xi=X_g$ and $\eta=X_H$.
Then we have that $f'\circ C=g+H\circ C$.
Since $f'\circ C(T,Q,Z,S,P)=f\circ \pi_{T,Q,P}\circ C^{-1}(t,q,z,a(T,Q,Z,0,P),p)=
f\circ \pi_{T,Q,P}(T,Q,Z,0,P)=f(T,Q,P)$ and 
$H$ has the form 
$H(t,q,z,s,p)=\sum_{i=1}^n h_i(t,q,z)p_i+h_0(t,q,z)+h'(t)s$,
We have that 
\[ f\equiv \sum_{i=1}^n (h_i(\Pi \circ C'))( p_i\circ C')+h_0(\Pi \circ C')+(h'(t\circ C'))(s\circ C')
\mbox{ mod } B_0.
 \]
Since $t\circ C=t$, we have the required form.\\
(a)$\Rightarrow$(u):
By Lemma \ref{g.function of rLu:lem}, there exists a 
function germ $H(T,Q,p)\in {\mathfrak M}(1+n+n)^2$
such that the function germ
$H(T,Q,p)-Ts$ is a generating function of $P_C$.
Then the function germ
$F(x,y,t,q,z)\in {\mathfrak M}(r;n+1+n+1)$ given by
$F(x,y,t,q,z)=-z+H(t,x,0,y)+\langle y,q\rangle$ is a generating family of ${\cal L}$.
Then $P'_C:=P_C|_{Z=S=0}$ has the form
\[ P'_C=\{ (T,Q,-\frac{\partial H}{\partial Q},
T,-\frac{\partial H}{\partial p},H-\langle 
\frac{\partial H}{\partial p},p\rangle ,\frac{\partial H}{\partial T},p)\}. \]
Then the map germ $P'_C\rightarrow ({\mathbb R}^{1+n+n},0),\ 
w\mapsto \pi_{T,Q,P}(w)$ is a diffeomorphism.
We set 
$D(F)=\{ (x,y,t,q,z)\in ({\mathbb H}^r\times
{\mathbb R}^{n+1+n+1},0)|F=x\frac{\partial F}{\partial x}=
\frac{\partial F}{\partial y}=0\}$.
We also define the map germ
$D(F)\rightarrow P'_C$ by 
\[
(x,y,t,q,z)\mapsto
(t,x,0,-\frac{\partial F}{\partial x},-\frac{\partial H}{\partial Q_{r+1}}(t,x,0,y),-\frac{\partial H}{\partial Q_n},t,q,z,\frac{\partial F}{\partial t},y).
\]
Then the composition of the above two map germs induces the 
map germ
${\cal E}_{T,Q,P}/B_0\rightarrow {\cal E}_{D(F)}$.
We denote $T,Q,P$ for the variables on the
source space of this map germ. 
Then the correspondence is given that:
\[
T\mapsto t,Q_1\mapsto x_1,\ldots,Q_r\mapsto x_r,
P_1\mapsto -\frac{\partial F}{\partial x_1},\ldots,
P_r\mapsto -\frac{\partial F}{\partial x_r},\]
\[ t\circ C'(T,Q,P)\mapsto t,
 q\circ C'(T,Q,P)\mapsto q,
z \circ C'(T,Q,P)\mapsto z,\]
\[ s\circ C'(T,Q,P)\mapsto \frac{\partial F}{\partial t},
p\circ C'(T,Q,P)\mapsto  y, ((\Pi\circ C')^*{\cal E}_{t,q,z})\mapsto {\cal E}_{t,q,z},
((t\circ C')^*{\cal E}_{t})\mapsto {\cal E}_{t}.
\]
Then (a) is transferred that  
\begin{eqnarray*}
{\cal E}(r;n+1+n+1)=
\langle F,x\frac{\partial F}{\partial x},
\frac{\partial F}{\partial y}\rangle_{{\cal E}(r;n+1+n+1)}\hspace{2cm}\\
+
\langle 1(=-\frac{\partial F}{\partial z}),y_1(=\frac{\partial F}{\partial q_1})
,\ldots,y_n(=\frac{\partial F}{\partial q_n})
\rangle_{{\cal E}_{t,q,z}}+
\langle \frac{\partial F}{\partial t} \rangle_{{\cal E}_{t}}. 
\end{eqnarray*}
It follows that $F$ is a reticular $t$-${\cal P}$-${\cal K}$-infinitesimal
stable unfolding of $F|_{t=0}$.
We have that $F$ is a reticular $t$-${\cal P}$-${\cal K}$-stable unfolding of  $F|_{t=0}$
by Theorem \ref{mthft:th}.
\hfill $\blacksquare$
\section{Genericity}\label{generic:sec}
\quad  
In order to give the generic classification of reticular Legendrian 
unfoldings,
we reduce our investigation to  finite dimensional jet-spaces 
of ${\cal P}$-contact diffeomorphism germs.
\begin{defn}\label{l,l+1detLeg}{\rm
Let ${\cal L}:({\mathbb L},0)\rightarrow (J^1({\mathbb R}\times 
{\mathbb R}^n,{\mathbb R}),0)$ be a reticular Legendrian unfolding.
We say that ${\cal L}$ is $l$-determined if the following condition holds:
For any extension $C\in C_T(J^1({\mathbb R}\times{\mathbb R}^n,{\mathbb R}),0)$ of ${\cal L}$, the reticular Legendrian unfolding 
$C'|_{{\mathbb L}}$ and ${\cal L}$ are ${\cal P}$-Legendrian 
equivalent for any $C'\in C_T(J^1({\mathbb R}\times{\mathbb R}^n,{\mathbb R}),0)$ satisfying that  $j^lC(0)=j^lC'(0)$.
}\end{defn}

As Lemma \ref{sta:lm}, we may consider the following other definitions
of finite determinacies of reticular Legendrian maps:\\
(1) The definition given by replacing $C_T(J^1({\mathbb R}\times{\mathbb R}^n,{\mathbb R}),0)$
to $C^\Theta_T(J^1({\mathbb R}\times{\mathbb R}^n,{\mathbb R}),0)$.\\
(2) The definition given by replacing $C_T(J^1({\mathbb R}\times{\mathbb R}^n,{\mathbb R}),0)$
to $C^Z_T(J^1({\mathbb R}\times{\mathbb R}^n,{\mathbb R}),0)$.\\
(3) The definition given by replacing $C_T(J^1({\mathbb R}\times{\mathbb R}^n,{\mathbb R}),0)$
to $C^{\Theta,Z}_T(J^1({\mathbb R}\times{\mathbb R}^n,{\mathbb R}),0)$.\\
Then the following holds:
\begin{prop}
Let ${\cal L}:({\mathbb L},0) \rightarrow (J^1({\mathbb R}\times{\mathbb R}^n,{\mathbb R}),0)$ be a reticular Legendrian unfolding.
Then \\
{\rm (A)} If ${\cal L}$ is $l$-determined of the original definition, then ${\cal L}$ is $l$-determined of the 
definition {\rm (1)}.\\
{\rm (B)} If ${\cal L}$ is $l$-determined of the definition {\rm (1)}, then ${\cal L}$ is $l$-determined of the definition {\rm (3)}.\\
{\rm (C)} If ${\cal L}$ is $l$-determined of the definition {\rm (3)}, then ${\cal L}$ is $(l+1)$-determined of the definition {\rm (2)}.\\
{\rm (D)} If ${\cal L}$ is $l$-determined of the definition {\rm (2)}, then 
${\cal L}$ is $l$-determined of the original definition.
\end{prop}
{\em Proof.} We need only to prove (C). 
Let $C\in C^Z_T(J^1({\mathbb R}\times{\mathbb R}^n,{\mathbb R}),0)$ be an 
extension of ${\cal L}$. 
Let $C'\in C^Z_T(J^1({\mathbb R}\times{\mathbb R}^n,{\mathbb R}),0)$ satisfying 
$j^{l+1}C(0)=j^{l+1}C'(0)$ be given.
Then there exist function germs $f(T,Q,S,P),g(T,Q,S,P)\in {\cal E}(2n+2)$
such that $C^{*}(dz-sdt-pdq)=-f(SdT+PdQ),C^{'*}(dz-sdt-pdq)=-g(SdT+PdQ)$.
Indeed $f$ is defined by that 
$fP_i=-\frac{\partial z_C}{\partial Q_i}+p_C
\frac{\partial q_C}{\partial Q_i}$ for $i=1,\ldots,n$, and 
$fS=-\frac{\partial z_C}{\partial T}+p_C
\frac{\partial q_C}{\partial T}$.
We define the diffeomorphism germs $\phi,\psi$ on 
$(J^1({\mathbb R}^n,{\mathbb R})\cap \{Z=0\},0)$ by
$\phi(T,Q,S,P)=(T,Q,fS,fP),\psi(T,Q,S,P)=(T,Q,gS,gP)$.
We set $C_1:=C\circ \phi^{-1},C'_1:=C'\circ\psi^{-1}\in 
C^{\Theta,Z}_T(J^1({\mathbb R}\times{\mathbb R}^n,{\mathbb R}),0)$
Then $j^l\phi(0)$ and $j^l\psi(0)$ depend only on $j^{l+1}C(0)$ and $j^{l+1}C'(0)$
respectively.
Therefore we have that $j^lC_1(0)=j^lC_1'(0)$.
Since ${\cal L}$ and $C_1|_{{\mathbb L}}$ are ${\cal P}$-Legendrian 
equivalent,
it follows that 
$C_1|_{{\mathbb L}}$ and $C'_1|_{{\mathbb L}}$ are  
${\cal P}$-Legendrian 
equivalent.
Therefore we have that ${\cal L}$ and $C'|_{{\mathbb L}}$ are 
${\cal P}$-Legendrian 
equivalent.\hfill $\blacksquare$\\

\begin{thm}\label{n+1det:Leg}
Let ${\cal L}:({\mathbb L},0) \rightarrow 
(J^1({\mathbb R}\times{\mathbb R}^n,{\mathbb R}),0)$ be a 
reticular Legendrian unfolding.
If ${\cal L}$ is infinitesimally stable then ${\cal L}$ is $(n+4)$-determined.
\end{thm}
{\em Proof.}
It is enough to prove ${\cal L}$ is $(n+3)$-determined of Definition 
\ref{l,l+1detLeg} (3).
Let $C\in C^{\Theta,Z}_T(J^1({\mathbb R}\times{\mathbb R}^n,{\mathbb R}),0)$ be an 
extension of ${\cal L}$.
Since the finite determinacy of reticular Legendrian unfoldings is
invariant under ${\cal P}$-Legendrian equivalences,
we may assume that $P_C$ has the form 
\[P_C=\{(T,Q,0,-\frac{\partial H}{\partial T}(T,Q,p)+s,
-\frac{\partial H}{\partial Q},T,-\frac{\partial H}{\partial p},
H-\langle \frac{\partial H}{\partial p},p\rangle ,s,p)\} \]
for some function germ $H(T,Q,p)\in {\mathfrak M}(2n+1)^2$.
Then the function germ $F(x,y,t,q,z)=-z+H_0(x,y,t)+\langle y,q\rangle  \in {\mathfrak M}(r;n+1+n+1)$
is a generating family of ${\cal L}$,
where $H_0(x,y,t)=H(t,x,0,y)\in {\mathfrak M}(r;n+1)^2$.
By Theorem \ref{staLeg:th}, we have that 
 $F$ is a
reticular $t$-${\cal P}$-${\cal K}$-stable unfolding of 
$f(x,y,q,z):=-z+H_0(x,y,0)+\langle y,q\rangle \in {\mathfrak M}(r;n+n+1)$.
Then $F$ is a
reticular $t$-${\cal P}$-${\cal K}$-infinitesimally stable unfolding of $f$ 
by Theorem \ref{mthft:th}.
This means that 
\[{\cal E}(r;n+1+n+1) =\langle F,x\frac{\partial
 F}{\partial x},\frac{\partial
 F}{\partial y}\rangle_{ {\cal
 E}(r;n+1+n+1) }
+\langle 1,\frac{\partial F}{\partial q}\rangle_{{\cal E}(1+n+1)}+
\langle 
\frac{\partial F}{\partial t}\rangle_{{\cal E}(1)}. 
\]
By the restriction of the both side to $q=z=0$, we have that
\begin{equation}
{\cal E}(r;n+1)=\langle H_0,x\frac{\partial
 H_0}{\partial x},\frac{\partial
 H_0}{\partial y}\rangle_{ {\cal
 E}(r;n+1) }
+\langle 1,y_1,\ldots,y_n,\frac{\partial
 H_0}{\partial t}\rangle_{{\cal E}(1)}.\label{rest:eqn}
\end{equation}
This means  that 
\begin{equation}
{\mathfrak M}(r;n+1)^{n+2}\subset \langle H_0,x\frac{\partial H_0}{\partial x},
\frac{\partial H_0}{\partial y}
\rangle_{{\cal E}(r;n+1)}+{\mathfrak M}(1){\cal E}(r;n+1).
\label{n+3det:leg}
\end{equation}
Let $C'\in C^{\Theta,Z}_T(J^1({\mathbb R}\times{\mathbb R}^n,{\mathbb R}),0)$ 
satisfying $j^{n+3}C(0)=j^{n+3}C'(0)$ be given.
There exists a function germ $H'(T,Q,p)\in {\mathfrak M}(2n+1)$ such that 
\[ P_{C'}=\{(T,Q,0,-\frac{\partial H'}{\partial T}(T,Q,p)+s,
-\frac{\partial H'}{\partial Q},T,-\frac{\partial H'}{\partial p},
H'-\langle \frac{\partial H'}{\partial p},p\rangle ,s,p)\}.\]
Since $H=z-qp$ on $P_C$ and $H'=z-qp$ on $P_{C'}$,
we have that $j^{n+3}H_0(0)=j^{n+3}H_0'(0)$, where 
$H_0'(x,y,t)=H'(t,x,0,y)\in {\mathfrak M}(r;n+1)^2$.
By (\ref{n+3det:leg}) we have that 
\[ {\mathfrak M}(r;n)^{n+2}\subset \langle H_0,x\frac{\partial H_0}{\partial x}
(x,y,0),
\frac{\partial H_0}{\partial y}(x,y,0)
\rangle_{{\cal E}(r;n)} \]
and this means that $H_0(x,y,0)$ is reticular ${\cal K}$-$(n+3)$-determined
by \cite[p.180 Lemma 2.3]{tPKfunct}.
Therefore we may assume that 
$H_0|_{t=0}=H'_0|_{t=0}$. It follows that
$H_0-H_0'\in {\mathfrak M}(1){\mathfrak M}(r;n+1)^{n+3}$.
Then the function germ $G(x,y,t,q,z)=-z+H_0'(x,y,t)+
\langle y,q\rangle \in {\mathfrak M}(r;n+1+n+1)$
is a generating family of $C'|_{{\mathbb L}}$.

We define the function germ $E_{\tau_0}(x,y,t,\tau)\in {\cal E}(r;n+1+1)$ by
$E_{\tau_0}(x,y,t,\tau)=(1-\tau-\tau_0)H_0(x,y,t)+(\tau+\tau_0)H_0'(x,y,t)$
for $\tau_0\in [0,1]$.
By (\ref{rest:eqn}) and (\ref{n+3det:leg}), we have that
\begin{equation}
{\mathfrak M}(r;n+1)^{n+3}
\subset \langle H_0,x\frac{\partial H_0}{\partial x}
\rangle_{{\cal E}(r;n+1)}+{\mathfrak M}(r;n+1)\langle 
\frac{\partial H_0}{\partial y}\rangle 
+{\mathfrak M}(1)\langle 1,y,\frac{\partial H_0}{\partial t}\rangle.\label{n+3eqn:eqn}
\end{equation}
Then we have that 
\begin{eqnarray*}
& & {\mathfrak M}_t{\mathfrak M}_{x,y,t}^{n+3}{\cal E}_{x,y,t,\tau} \\
& = & {\mathfrak M}_t{\mathfrak M}_{x,y,t}^{n+3}
({\cal E}_{x,y,t}+{\mathfrak M}_\tau{\cal E}_{x,y,t,\tau})\\
& \subset & 
{\mathfrak M}_t{\mathfrak M}_{x,y,t}^{n+3}
+{\mathfrak M}_t{\mathfrak M}_\tau{\mathfrak M}_{x,y,t}^{n+3}{\cal E}_{x,y,t,\tau}\\
& \subset & 
{\mathfrak M}_t(
\langle H_0,x\frac{\partial H_0}{\partial x}\rangle_{{\cal E}_{x,y,t}}+
{\mathfrak M}_{x,y,t}\langle\frac{\partial H_0}{\partial y}
\rangle +{\mathfrak M}_t \langle 1,y,\frac{\partial
 H_0}{\partial t}\rangle )
 +{\mathfrak M}_{t,\tau}{\mathfrak M}_t
{\mathfrak M}_{x,y,t}^{n+3}{\cal E}_{x,y,t,\tau}\\
& \subset & 
{\mathfrak M}_{t,\tau}
\langle E_{\tau_0},x\frac{\partial E_{\tau_0}}{\partial x}\rangle_{{\cal E}_{x,y,t,\tau}}+{\mathfrak M}_{t,\tau}{\mathfrak M}_{x,y,t,\tau}\langle\frac{\partial E_{\tau_0}}{\partial y}
\rangle \\
& & \hspace{5cm}+{\mathfrak M}_{t,\tau}^2 \langle 1,y,\frac{\partial
 E_{\tau_0}}{\partial t}\rangle +{\mathfrak M}_{t,\tau}{\mathfrak M}_t
{\mathfrak M}_{x,y,t}^{n+3}{\cal E}_{x,y,t,\tau}.
\end{eqnarray*}
By (\ref{n+3eqn:eqn}) we have the second inclusion.
For the last inclusion, observe that
\[ x_j\frac{\partial E_{\tau_0}}{\partial x_j}-x_j\frac{\partial H_0}{\partial x_j}=
(\tau_0+\tau)x_j\frac{\partial}{\partial x_j}(H_0'-H_0)\in {\mathfrak M}_t
{\mathfrak M}_{x,y,t}^{n+3}, \]
\[ \frac{\partial E_{\tau_0}}{\partial y_j}-\frac{\partial H_0}{\partial y_j}=
(\tau_0+\tau)\frac{\partial}{\partial y_j}(H_0'-H_0)\in {\mathfrak M}_t
{\mathfrak M}_{x,y,t}^{n+2}, \]
\[ \frac{\partial E_{\tau_0}}{\partial t}-\frac{\partial H_0}{\partial t}=
(\tau_0+\tau)\frac{\partial}{\partial t}(H_0'-H_0)\in
{\mathfrak M}_{x,y,t}^{n+3}. \]
By Malgrange preparation theorem we have that 
\begin{eqnarray*}
\frac{\partial E_{\tau_0}}{\partial \tau}\in {\mathfrak M}_t{\mathfrak M}_{x,y,t}^{n+3}\subset 
{\mathfrak M}_t{\mathfrak M}_{x,y,t}^{n+3}{\cal E}_{x,y,t,\tau}\hspace{5cm}
\\
\subset  {\mathfrak M}_{t,\tau}(
\langle E_{\tau_0},x\frac{\partial E_{\tau_0}}{\partial x}\rangle_{{\cal E}_{x,y,t,\tau}}+
{\mathfrak M}_{x,y,t,\tau}\langle\frac{\partial E_{\tau_0}}{\partial y}
\rangle) 
+{\mathfrak M}_{t,\tau}^2 \langle 1,y,\frac{\partial
 E_{\tau_0}}{\partial t}\rangle.
\end{eqnarray*}
for $\tau_0\in [0,1]$.
By using analogous methods of \cite[p.37 Theorem 2.6]{staun},
we have that there exist $\Phi(x,y,t)\in {\cal B}_1(r;n+1)$, a unit $a\in
{\cal E}(r;n+1)$ and $b_1(t),\ldots,b_n(t),c(t)\in {\mathfrak M}(1)$ such that\\
(1) $\Phi$ has the form:
$
\Phi(x,y,t)=(x\phi_1(x,y,t),\phi_2(x,y,t),\phi_3(t))
$,\\
(2) $H_0(x,y,t)=a(x,y,t)\cdot H_0'\circ \Phi(x,y,t)+\sum_{i=1}^ny_ib_i(t)+c(t)$ 
for $(x,y,t)\in ({\mathbb H}^r\times {\mathbb R}^{n+1},0)$

We define the reticular $t$-${\cal P}$-${\cal K}$-isomorphism $(\Psi,d)$
by 
\[
\Psi(x,y,t,q,z)=(x\phi_1(x,y,t),\phi_2(x,y,t),\phi_3(t),q(1-b(t)),z),
d(x,y,t,q,z)=a(x,y,t).
\]
We set $G':=d\cdot G\circ\Psi\in {\mathfrak M}(r;n+n+1)$.
Since $\frac{\partial E_{\tau_0}}{\partial \tau}|_{t=0}=0$, we have that 
$a(x,y,0)=1$ and $\Phi(x,y,0)=(x,y,0)$.
Therefore we have that $G'|_{t=0}=f$.
Then $F$ and $G'$ are reticular $t$-${\cal P}$-${\cal K}$-infinitesimal versal
unfoldings of $F|_{t=0}$.
Since $G$ and $G'$ are reticular $t$-${\cal P}$-${\cal K}$-equivalent,
it follows that $F$ and $G$ are reticular $t$-${\cal P}$-${\cal K}$-equivalent.
Therefore ${\cal L}$ and $C'|_{{\mathbb L}}$ are
${\cal P}$-Legendrian equivalent.\hfill $\blacksquare$\\

Let  $J^l(2n+3,2n+3)$ be the set of $l$-jets of map germs  from $(J^1({\mathbb R}\times{\mathbb R}^n,{\mathbb R}),0)$ 
to $(J^1({\mathbb R}\times{\mathbb R}^n,{\mathbb R}),0)$ and 
$pC^l(n)$ be the Lie group in $J^l(2n+3,2n+3)$ which
consists of $l$-jets of ${\cal P}$-contact diffeomorphism germs 
on $(J^1({\mathbb R}\times{\mathbb R}^n,{\mathbb R}),0)$.
Let $L^l(2n+3)$ be the Lie group which consists of $l$-jet of diffeomorphism germs on 
$(J^1({\mathbb R}\times{\mathbb R}^n,{\mathbb R}),0)$.

We  consider the Lie subgroup $rpLe^l(n)$ of $L^l(2n+3)\times L^l(2n+3)$ which consists of 
 $l$-jets of reticular ${\cal P}$-diffeomorphisms on the  source space and $l$-jets of ${\cal P}$-Legendrian 
equivalences of $\Pi$ on the target space: 
\begin{eqnarray*}
 rpLe^l(n)=\{ (j^l\phi(0),j^lK(0))\in L^l(2n+3)\times L^l(2n+3)\ |\ \phi
 \mbox{ is a reticular} \hspace{2cm}\\
\mbox{ ${\cal P}$-diffeomorphism on }
(J^1({\mathbb R}\times{\mathbb R}^n,{\mathbb R}),0),  
K \mbox{ is a ${\cal P}$-Legendrian equivalence of } 
\Pi \}.
\end{eqnarray*} 
The group $rpLe^l(n) $ acts on $J^l(2n+3,2n+3)$ and $pC^l(n)$ is invariant 
under this action.
Let $C$ be  a ${\cal P}$-contact diffeomorphism germ on 
$(J^1({\mathbb R}\times{\mathbb R}^n,{\mathbb R}),0)$ and set 
$z=j^lC(0)$, ${\cal L}=C|_{{\mathbb L}}$.
We denote  the orbit $rpLe^l(n)\cdot z$ by $[z]$.
Then 
\[ [z]=\{ j^lC'(0)\in pC^l(n)\ | \ {\cal L}  \mbox{ and } C'|_{{\mathbb L}} 
\mbox{ are ${\cal P}$-Legendrian equivalent} \}. \]

We denote by $VI_C$  the vector space which consists of infinitesimal ${\cal P}$-contact transformation germs of $C$
and  denote by $VI^0_C$ the subspace of $VI_C$ which consists of germs which
 vanish on $0$.
We denote by $VL_{J^1({\mathbb R}\times{\mathbb R}^n,{\mathbb R})}$  the vector space consisting of 
infinitesimal ${\cal P}$-Legendrian 
equivalences of $\Pi$  and denote by $VL^0_{J^1({\mathbb R}\times{\mathbb R}^n,{\mathbb R})}$
the subspace of $VL_{J^1({\mathbb R}\times{\mathbb R}^n,{\mathbb R})}$ consists of germs which vanish at $0$.
We denote by $V^0_{{\mathbb L}}$ the vector space consists of infinitesimal reticular ${\cal P}$-diffeomorphisms 
on $(J^1({\mathbb R}\times{\mathbb R}^n,{\mathbb R}),0)$ which vanishes at $0$.
By Lemma \ref{infsta:t-Leglem}, we have that:
\begin{eqnarray*}
VI^0_C=\{ v:(J^1({\mathbb R}\times{\mathbb R}^n,{\mathbb R}),0)\rightarrow (T(J^1({\mathbb R}\times{\mathbb R}^n,{\mathbb R})),0)\ | \hspace{2cm}\\
v=X_f\circ C \mbox{ for some }f\in {\mathfrak M}^2_{t,q,z,p}\},\\ 
VL^0_{J^1({\mathbb R}\times{\mathbb R}^n,{\mathbb R})}=\{\eta\in X(J^1({\mathbb R}\times{\mathbb R}^n,{\mathbb R}),0)\ | \hspace{5cm}\\
\eta=X_H
\mbox{ for some ${\cal P}$-fiber preserving function germ } H\in {\mathfrak M}^2_{J^1({\mathbb R}\times{\mathbb R}^n,{\mathbb R})}\},\\
V^0_{{\mathbb L}}=\{ \xi \in X(J^1({\mathbb R}\times{\mathbb R}^n,{\mathbb R}),0)\ | 
\xi=X_g  \mbox{ for some }g\in B'\},
\end{eqnarray*}
where $B'=\langle q_1p_1,\ldots,q_rp_r\rangle_{J^1({\mathbb R}\times{\mathbb R}^n,{\mathbb R})}+
{\mathfrak M}_{J^1({\mathbb R}\times{\mathbb R}^n,{\mathbb R})}
\langle 
q_{r+1},\ldots,q_n,z\rangle
+{\mathfrak M}_t\langle s \rangle$.\vspace{2mm}

We define the homomorphisms $tC:V^0_{{\mathbb L}}\rightarrow VI^0_C$ by 
$tC(v)=C_*v$ and 
$wC:VL^0_{J^1({\mathbb R}\times{\mathbb R}^n,{\mathbb R})}\rightarrow 
VI^0_C$ by 
$wC(\eta)=\eta\circ C$.

\vspace{2mm}
We denote $VI^{l}_C$ the subspace of $VI_C$ consists of infinitesimal
${\cal P}$-contact transformation germs of $C$ whose $l$-jets are $0$:
\[  VI^l_C=\{ v\in VI_C\ | \ j^lv(0)=0\}. \]

For $\tilde{C}\in C^\Theta_T(U,J^1({\mathbb R}\times{\mathbb R}^n,{\mathbb R}))$,
we define
 the continuous map 
$j^l_0\tilde{C}:U\rightarrow pC^l(n)$ by mapping $w=(T^0,Q^0,Z^0,S^0,P^0)$ to the $l$-jet of the 
${\cal P}$-contact diffeomorphism germ $\tilde{C}_w$ at $0$, 
where $\tilde{C}_w(x)$ is given by $x=(T,Q,Z,S,P)\mapsto \tilde{C}(w+x)-\tilde{C}(w)+(0,0,s_{\tilde{C}}(w)t_{\tilde{C}}(w+x)-s_{\tilde{C}}(w)t_{\tilde{C}}(w)+
p_{\tilde{C}}(w)q_{\tilde{C}}(w+x)-p_{\tilde{C}}(w)q_{\tilde{C}}(w)+S^0T+P^0Q,0,0)$.

We also define $j^l_0C:(J^1({\mathbb R}\times{\mathbb R}^n,{\mathbb R}),0)\rightarrow pC^l(n)$ by the same method for $C\in C_T^\Theta(J^1({\mathbb R}\times{\mathbb R}^n,{\mathbb R}),0)$.
\begin{prop}
Let  $C\in C_T^\Theta(J^1({\mathbb R}\times{\mathbb R}^n,{\mathbb R}),0)$ and set 
$z=j^lC(0)$.
Then  $j^l_0C$ is transversal to $[z]$ if and only if
\begin{equation}
tC( V^0_{{\mathbb L}})+wC(VL_{J^1({\mathbb R}\times{\mathbb R}^n,{\mathbb R})}))+VI^{l+1}_C=VI_C. 
\label{trans:eqn}
\end{equation}
\end{prop}
{\em Proof}.
We consider the surjective projection $\pi_l:VI_C\rightarrow T_z(pC^l(n))$. 
Since $(j^lC)_*(v)=\pi_l( C_*v)$ for all $v\in T_0(J^1({\mathbb R}\times{\mathbb R}^n,{\mathbb R}))$,
it follows that 
$j^lC$ is transversal to $[z]$ if and only if 
$(j^lC)_*(T_0(J^1({\mathbb R}\times{\mathbb R}^n,{\mathbb R})))+T_z[z]=T_z(pC^l(n))$ and this holds 
if and only if 
\[
 (\pi_l)^{-1}((j^lC)_*(T_0(J^1({\mathbb R}\times{\mathbb R}^n,{\mathbb R})))+
tC( V^0_{{\mathbb L}})
+wC(VL^0_{J^1({\mathbb R}\times{\mathbb R}^n,{\mathbb R})}))=VI_C
\]
 and 
this holds if and only if (\ref{trans:eqn}) holds.
\hfill $\blacksquare$\\
\begin{thm}\label{stabletrans_tleg:th}
Let ${\cal L}$ be a reticular Legendrian unfolding. 
Let $C\in C^\Theta_T(J^1({\mathbb R}\times{\mathbb R}^n,{\mathbb R}),0)$ be an extension of ${\cal L}$ and $l\geq (n+2)^2$.
Then the followings are equivalent:\\
{\rm (s)} ${\cal L}$  is stable.\\
{\rm (t)} $j^l_0C$ is transversal to $[j^l_0C(0)]$.\\
{\rm (a')} ${\cal E}_{t,q,p}=
B_0+
\langle 1,p_1\circ C',\ldots,p_n\circ C'\rangle_{(\Pi\circ C')^*{\cal E}_{t,q,z}}+
\langle s\circ C'\rangle_{{\cal E}_t}+{\mathfrak M}_{t,q,p}^l$,
where $C'=C|_{z=s=0}$ and $B_0=\langle q_1p_1,\ldots,q_rp_r,
q_{r+1},\ldots,q_n\rangle_{{\cal E}_{t,q,p}}$.\\
{\rm (a)} ${\cal E}_{t,q,p}=
B_0+
\langle 1,p_1\circ C',\ldots,p_n\circ C'\rangle_{(\Pi\circ C')^*{\cal E}_{t,q,z}}+
\langle s\circ C'\rangle_{{\cal E}_t}$.\\
{\rm (is)} ${\cal L}$ is infinitesimally stable. 
\end{thm}
{\em Proof}.
(s)$\Rightarrow$(t): Let $\tilde{C}\in 
C^\Theta_T(U,J^1({\mathbb R}\times{\mathbb R}^n,{\mathbb R}))$ be a representative 
of $C$.
By theorem \ref{t_contra:th} and (s), there exists  $\tilde{C'}
\in C^\Theta_T(U,J^1({\mathbb R}\times{\mathbb R}^n,{\mathbb R}))$ such that 
$j^l_0{\tilde{C'}}$ is transversal to $[j^l_0C(0)]$, and
$\tilde{C'}_w|_{{\mathbb L}}$ and ${\cal L}$ are
${\cal P}$-Legendrian equivalent for some $w\in U$.
This means that $[j^l_0\tilde{C'_w}(0)]=[j^lC(0)]$ and hence 
$j^l_0C$ is transversal to $[j^l_0C(0)]$ at $0$.\\
(t)$\Leftrightarrow$(a): This is proved by 
an analogous method of Theorem \ref{staLeg:th}.\\
(a)$\Leftrightarrow$(a'): We need only to prove (a')$\Rightarrow$(a).
By the restriction of  (a') to $t=0$ we have that:
\[{\cal E}_{q,p}=
B'_0+
\langle 1,p_1\circ C'',\ldots,p_n\circ C''\rangle_{(\Pi\circ C'')^*{\cal E}_{t,q,z}}+
\langle s\circ C''\rangle_{{\mathbb R}}+
{\mathfrak M}_{q,p}^l,
\]
where $C''=C'|_{t=0}$ and $B'_0=B_0|_{t=0}$.
Then we have that 
\[ {\cal E}_{q,p}=B'_0+(\Pi\circ C'')^*
{\mathfrak M}_{t,q,z}{\cal E}_{q,p}+
\langle 1,p_1\circ C'',\ldots,p_n\circ C'',
s\circ C''\rangle_{{\mathbb R}}+{\mathfrak M}_{q,p}^l.\]
It follows that 
\[ {\mathfrak M}_{q,p}^{n+2}\subset B'_0+(\Pi\circ C'')^*
{\mathfrak M}_{t,q,z}{\cal E}_{q,p}.\]
Therefore 
\[ {\mathfrak M}_{t,q,p}^{n+2}\subset B_0+(\Pi\circ C')^*
{\mathfrak M}_{t,q,z}{\cal E}_{q,p}+
{\mathfrak M}_{t}{\cal E}_{t,q,p},\]
and we have that 
\[ {\mathfrak M}_{t,q,p}^l=({\mathfrak M}_{t,q,p}^{n+2})^{n+2}
\subset B_0+(\Pi\circ C')^*
{\mathfrak M}_{t,q,z}^{n+2}{\cal E}_{t,q,p}+
{\mathfrak M}_{t}{\cal E}_{t,q,p}.\]
It follows that
\[
{\cal E}_{t,q,p}=
B_0+
\langle 1,p_1\circ C',\ldots,p_n\circ C'\rangle_{(\Pi\circ C')^*{\cal E}_{t,q,z}}+
\langle s\circ C'\rangle_{{\cal E}_t}+(\Pi\circ C')^*
{\mathfrak M}_{t,q,z}^{n+2}{\cal E}_{t,q,p}+
{\mathfrak M}_{t}{\cal E}_{t,q,p}.
\]
This means (a) by \cite[Corollary 1.8]{spsing}.\\
(a)$\Leftrightarrow$(is): This is proved in Theorem \ref{staLeg:th}.\\
(t)\&(is)$\Rightarrow$(s):
Since $j^l_0C$ is transversal to $[j^l_0C(0)]$, it follows that 
there exist
a representative $\tilde{C}\in C^\Theta_T(U,J^1({\mathbb R}\times{\mathbb R}^n,{\mathbb R}))$ of $C$ and
a neighborhood  $W_{\tilde{C}}$ of $\tilde{C}$ in 
$C^\Theta_T(U,J^1({\mathbb R}\times{\mathbb R}^n,{\mathbb R}))$ 
such that for any $\tilde{C'}\in W_{\tilde{C}}$
 there exists $w\in U$ such that 
$j^l_0\tilde{C}'$ is transversal to  $[j^l_0C(0)]$ at $w$.
Since $j^l_0\tilde{C'_w}(0)\in [j^l_0C(0)]$, it follows that there exists 
$C''\in C_T^\Theta(J^1({\mathbb R}\times{\mathbb R}^n,{\mathbb R}),0)$ such that ${\cal L}$ and 
$C''|_{{\mathbb L}}$ are 
${\cal P}$-Legendrian equivalent
and $j^l_0C''(0)=j^l_0\tilde{C'_w}(0)$. 
Since ${\cal L}$ is infinitesimally stable, it follows that ${\cal L}$ is $(n+4)$-determined
by Theorem \ref{n+1det:Leg}.
Therefore we have that $C''|_{{\mathbb L}}$ is also 
$(n+4)$-determined.
Then $C''|_{{\mathbb L}}$ and
 $\tilde{C'_w}|_{{\mathbb L}}$ are ${\cal P}$-Legendrian equivalent.
This means that ${\cal L}$ is stable.
\hfill $\blacksquare$
%
%
\vspace{2mm}

Let ${\cal L}$ be a stable reticular Legendrian unfolding.
We say that ${\cal L}$ is {\em simple} if 
there exists a representative $\tilde{C}\in C_T^\Theta(
U,J^1({\mathbb R}\times{\mathbb R}^n,{\mathbb R}))$ of
 a extension of ${\cal L}$ such that 
$\{ \tilde{C}_w| w\in  U\}$ is covered by 
finite orbits $[C_1],\ldots,[C_m]$ for  
 some $C_1,\ldots,C_m\in C_T(J^1({\mathbb R}\times{\mathbb R}^n,{\mathbb R}),0)$. 

\begin{lem}\label{simpleLegket:lem}
Let ${\cal L}$
be a stable reticular Legendrian unfolding
and $l\geq (n+2)^2$.
Let $C\in C_T^\Theta(J^1({\mathbb R}\times{\mathbb R}^n,{\mathbb R}),0)$ be an extension of $\cal{L}$.
Then ${\cal L}$ is simple if and only if there exists an open neighborhood $W_z$ of 
$z=j^{l}_0C(0)$ in $pC^{l}(n)$ and $z_1,\ldots,z_m\in pC^{l}(n)$ such 
that 
$W_z\subset [z_1]\cup\cdots\cup [z_m]$.
\end{lem}
{\em Proof}.
Suppose that ${\cal L}$ is simple.
Then there exists a representative 
$\tilde{C}\in C_T^\Theta(U,J^1({\mathbb R}\times{\mathbb R}^n,{\mathbb R}))$ 
of
 a extension of ${\cal L}$ and $C_1,\ldots,C_m\in C_T(J^1({\mathbb R}\times{\mathbb R}^n,{\mathbb R}),0)$
such that 
\begin{equation}
\{ \tilde{C}_w| w\in U\}\subset 
[C_1]\cup\cdots\cup [C_m].\label{trans1:eqn}
\end{equation}
Since ${\cal L}$ is stable, it follows that 
$j^{l}_0\tilde{C}$ is transversal to $[z]$ at $0$ by 
Theorem \ref{stabletrans_tleg:th}.
This means that 
there exists a neighborhood $W_z$ of $z$ in $pC^l(n)$ such that 
$W_z\subset \cup_{w\in U}[j^{l}_0\tilde{C}(w)]$.
It follows that  $W_z
\subset 
[j^{l}C_1(0)]\cup\cdots\cup [j^{l}C_m(0)]$.

Conversely suppose that 
there exist a neighborhood $W_z$ of $z$
 in $pC^{l}(n)$
and $z_1,\ldots,z_m
\in pC^{l}(n)$ such that 
$W_z\subset [z_1]\cup\cdots\cup [z_m]$.
Since the map $j^{l}_0\tilde{C}:U\rightarrow pC^{l}(n)$ is 
continuous,
there exists a neighborhood $U'$ of $0$ in $U$
such that  $j^{l}_0\tilde{C}(w)\in W_z$ for any $w\in {U'}$.
Then we have that 
$\cup_{w\in {U'}}j^{l}_0\tilde{C}(x)
\subset [z_1]\cup\cdots\cup [z_m]$.
Choose ${\cal P}$-contact  diffeomorphism germs $C_1,\ldots,C_m$ on 
$(J^1({\mathbb R}\times{\mathbb R}^n,{\mathbb R}),0)$ such that 
$j^{l}C_j(0)=z_i$ 
for $i=1,\ldots,m$.
By Theorem \ref{stabletrans_tleg:th} (a'), we may assume that 
 each reticular Legendrian unfolding 
$C_i|_{{\mathbb L}}$ is stable, thus
$l$-determined.
For any $w\in {U'}$ 
there exists $i\in \{1,\ldots,m\}$ such that
$j^{l}_0\tilde{C}(w)\in [j^{l}C_i(0)]$.
It follows that reticular Legendrian unfoldings 
$\tilde{C_w}|_{{\mathbb L}}$  and
 $C_i|_{{\mathbb L}}$
are ${\cal P}$-Legendrian equivalent.
Therefore $\tilde{C}_w\in [C_i]$.
We have (\ref{trans1:eqn}). \hfill $\blacksquare$\\

 \begin{lem}\label{simpleLeg:lem}
A stable reticular Legendrian  unfolding ${\cal L}$ is simple if and only if 
for a generating family 
$F(x,y,t,q,z)\in {\mathfrak M}(r;k+1+n+1)$ of ${\cal L}$,
$f(x,y)=F(x,y,0,0)\in {\mathfrak M}(r;k)^2$ is a  ${\cal K}$-simple singularity.
\end{lem}
{\em Proof.}
Suppose that ${\cal L}$ is simple. 
Then $\Pi\circ {\cal L}$ is simple as a reticular Legendrian map.
It follows that $f$ is a ${\cal K}$-simple singularity by 
\cite{generic}.

Conversely suppose that $f$ is ${\cal K}$-simple.
Let $l\geq (n+2)^2$.
There exist a neighborhood $W_z$ of $z=j^lf(0)$ in $J^l(r+k,1)$
and $f_1,\ldots,f_m\in {\cal E}(r;k)$ such that 
$W_z\subset [j^lf_1(0)]\cup\cdots\cup[j^lf_m(0)]$.
By the ${\cal K}$-simplicity of $f$, we can choose $f_i's$ such  that
each $f_i$ is simple and has reticular ${\cal K}$-codimension$\leq n+2$, and hence
is reticular ${\cal K}$-$l$-determined.
We choose a reticular  $t$-${\cal P}$-${\cal K}$-stable unfolding 
$F^1_i(x,y,t,q,z)\in {\mathfrak M}(r;k+1+n+1)$ of $f_i$ for each $i$.
If  there exists a reticular  ${\cal P}$-${\cal K}$-stable unfolding as of 
$f_i$ as $(n+2)$-dimensional unfolding, 
we set it by $F^0_i$, otherwise set $F^0_i=F^1_i$ for each $i$.
We also choose an extension $C^j_i\in C_T(J^1({\mathbb R}\times{\mathbb R}^n,{\mathbb R}),0)$ of a reticular Legendrian unfolding of which
$F^j_i$ is a generating family for each $i,j$.
Let $C_0$ be an extension of ${\cal L}$.
We may assume that the canonical relation $P_{C_0}$ has the form:
\[
P_{C_0}=\{(T,Q,Z,-\frac{\partial H_{C_0}}{\partial T}(T,Q,p)+s,
-\frac{\partial H_{C_0}}{\partial Q},T,
-\frac{\partial H_{C_0}}{\partial p},
H_{C_0}-\langle \frac{\partial H_{C_0}}{\partial p},p\rangle+Z,s,p)\}
\]
for some function germ $H_{C_0}$.
Then the $l$-jet of $H_{C_0}$ is determined by the $l$-jet of $C_0$ 
since $H_{C_0}=z-qp$ on $P_{C_0}$.
For a ${\cal P}$-contact diffeomorphism germ $C$ on 
$(J^1({\mathbb R}\times{\mathbb R}^n,{\mathbb R}),0)$ around $C_0$,
there exists a function germ $H_C(T,Q,p)$ satisfying the above condition
for $P_C$.
We define $H'_C(x,y)\in {\mathfrak M}(r;n)$ by
$H'_C(x,y)=H_C(0,x,0,y)$.
Then there exists a neighborhood $U$ of $j^lC_0(0)$ such that 
 the following continuous maps are constructed: 
\[ \begin{array}{ccccccc}
U & \rightarrow & J^l(1+n+n,1)& \rightarrow & J^l(r+n,1)\\
j^lC(0) & \mapsto & j^lH_C(0) & \mapsto & j^lH'_C(0)
\end{array}.\]
Since $H'_{C_0}$ is reticular ${\cal K}$-equivalent to $f$, 
we may assume that $j^lH'_{C_0}(0)\in W_z$.
We set $U'$ the inverse image of $W_z$ by the above maps.
For any $j^lC(0)\in U'$,
there exists a number $i$ such that $j^lH'_C(0)\in [f_i]$.
By Theorem \ref{stabletrans_tleg:th} (a'), we have that 
$C|_{{\mathbb L}}$ is a stable reticular Legendrian unfolding.
Since $f_i$ is reticular ${\cal K}$-$l$-determined,
we have that $H'_C$ and $f_i$ are reticular ${\cal K}$-equivalent.
Then the reticular Legendrian unfolding $C|_{{\mathbb L}}$ is 
${\cal P}$-Legendrian equivalent to $C^0_i|_{{\mathbb L}}$ or 
$C^1_i|_{{\mathbb L}}$ and it follows that 
$j^lC(0)\in [C^0_i]\cup [C^1_i]$.
Then we have that 
\[ U_z\subset [C^0_1]\cup [C^1_1]\cup\cdots \cup [C^0_m]\cup [C^1_m]\]
and this means that ${\cal L}$ is simple.\hfill $\blacksquare$\\

By \cite[Proposition 6.5]{tPKfunct}, we have that:
\begin{thm}\label{classfunct:tth}
Let $F(x,y,t,q,z)\in {\mathfrak M}(r;k+n+1)$ be a 
${\cal P}$-$C$-non-degenerate function germ for $r=0,n\leq 4$ or 
$r=1,n\leq 2$.
Then $F$ is stably reticular 
$t$-${\cal P}$-${\cal K}$-equivalent for one of the following types.\\
In the case $r=0,n\leq 4$:
$({}^0A_l)\ y_1^{l+1} +\displaystyle{\sum_{i=1}^{l}}
q_iy_1^i+z\ (2\leq l\leq n)$,\\
$({}^0D^\pm_4)\  y_1^2y_2\pm y_2^3+q_1y_2^2+q_2y_2+q_3y_1
+z,$\\
$({}^0D_5)\ y_1^2y_2+ y_2^4+q_1y_2^3+q_2y_2^2+q_3y_2+q_4y_1+z$,\\
$({}^1A_l)\ y_1^{l+1} +
(t+q_{l}^2\pm q_{l+1}^2\pm \cdots\pm q_n^2)y_1^{l-1}+
\displaystyle{\sum_{i=1}^{l-1}}
q_iy_1^i+z\ (3\leq l\leq n)$,\\
$({}^1D^\pm_4)\ y_1^2y_2\pm y_2^3+
ty_2^2+q_1y_2+q_2y_1+z,\ 
y_1^2y_2\pm y_2^3+
(t+ q_3^2)y_2^2+q_1y_2+q_2y_1+z$,\\
 $({}^1D_5)\
y_1^2y_2+ y_2^4+ty_2^3+q_1y_2^2+q_2y_2+q_3y_1+z,\
 y_1^2y_2+ y_2^4+(t+ q_4^2 )y_2^3+q_1y_2^2+q_2y_2+q_3y_1
+z$, \\
 $({}^1D^\pm_6)\ y_1^2y_2\pm  y_2^5+ty_2^6+q_1y_2^3+q_2y_2^2+q_3y_2+q_4y_1+z$,\\
$({}^1E_6)\ y_1^3+ y_2^4+ty_1y_2^2+q_1y_1y_2+q_2y_2^2+
q_3y_1+q_4y_2+z$.\\
In the case $r=1,n\leq 2$:
$({}^0B_2)\ x^2+q_1x+z$,\\
$({}^0B_3)\ x^3+q_1x^2+q_2x+z$,\\
$({}^0C^\pm_3)\ \pm xy+y^3+q_1y^2+q_2y+z$,\\
$({}^1B_3)\ x^3+tx^2+q_1x+z,\ 
x^3+(t\pm q_2^2)x^2+q_1x+z$,\\
$({}^1B_4)\ x^4+tx^3+q_1x^2+q_2x+z$,\\
$({}^1C^\pm_3)\ \pm xy+y^3+ty^2+
q_1y+z,\ 
\pm xy+y^3+(t+ q_2^2)y^2+
q_1y+z$,\\
$({}^1C_4)\ xy+y^4+ty^3+q_1y^2+q_2y+z$,\\
$({}^1F_4)\ x^2+y^3+txy+q_1x+q_2y+z$.\\
\end{thm}
\begin{thm}\label{genericclass:tth}
Let $r=0,n\leq 4$ or $r=1,n\leq 2$.
Let  $U$ be a neighborhood of $0$ in 
$J^1({\mathbb R}\times{\mathbb R}^n,{\mathbb R})$.
Then there exists a residual set $O\subset  C^\Theta_T(U,J^1({\mathbb R}\times{\mathbb R}^n,{\mathbb R}))$ 
such that for any $\tilde{C}\in O$ and $w\in U$,
the reticular Legendrian unfolding $\tilde{C}_w|_{{\mathbb L}}$ 
is stable and has a generating family which is stably reticular 
$t$-${\cal P}$-${\cal K}$-equivalent for one of the types in the 
previous theorem.
\end{thm}
{\em Proof}.
In the case $r=1,\ n\leq 2$.
Let $F_X(x,y,t,q)\in {\mathfrak M}(r;k+1+n)$ 
be a reticular $t$-${\cal P}$-${\cal K}$-stable unfolding of singularity 
$X\in {\mathfrak M}(r;k)^2$ for 
\[ X=B_2,B_3,B_4,
C^\pm_3,C_4,F_4.\]
Then other unfoldings are not stable since other singularities have 
reticular ${\cal K}$-codimension $>4$.
We choose stable reticular Legendrian unfoldings ${\cal L}_X:
({\mathbb L},0)\rightarrow 
(J^1({\mathbb R}\times {\mathbb R}^n,{\mathbb R}),0)$
with the generating family $F_X$,
 and $C_X$ be an extension of ${\cal L}_X$ for above list.
Let $l>16$. We define that 
\[
 O'=\{ \tilde{C}\in  
C^\Theta_T(U,J^1({\mathbb R}\times{\mathbb R}^n,{\mathbb R}))\ |
j^{l}_0\tilde{C} \mbox{ is transversal to }[j^{l}C_X(0)]
\mbox{ for all } X \}.
\]
Then $O'$ is a residual set. We set 
\[ Y=\{ j^{l}C(0)\in C^{l}(n)\ |\ 
\mbox{the codimension of }[j^{l}C(0)]>2n+4\}.\]
Then $Y$ is an algebraic set in $pC^{l}(n)$.
Therefore we can define that 
\[ O''=\{ \tilde{C}\in  C^\Theta_T(U,J^1({\mathbb R}\times{\mathbb R}^n,{\mathbb R}))\ |\ 
j^{l}_0\tilde{C}\mbox{ is transversal to } Y\}.\]
Then $Y$ has codimension $>2n+4$ because all 
${\cal P}$-contact diffeomorphism germ with 
$j^{l}C(0)\in Y$ adjoin to the above list which are simple.
Therefore the set
\[ O''=\{ \tilde{C}\in  C^\Theta_T(U,J^1({\mathbb R}\times{\mathbb R}^n,{\mathbb R}))\ |\ 
j^{l}_0\tilde{C}(U)\cap Y=\emptyset \}\]
is residual.
Then the set $O=O'\cap O''$ has the required condition.\hfill $\blacksquare$
\section{Classifications for the cases $r\geq 2$}\label{classover1}
\quad
In order to classify generic bifurcations of wavefronts for the case $r\geq2$, we can not use
our equivalence relation.
For example, consider the bifurcations with the generating families:
$F_a(x_1,x_2,t,q_1,q_2,q_3)=x_1^2+a x_1x_2+x_2^2+tx_1x_2+q_1x_1+q_2x_2+q_3(a^2\neq 4)$.
Then $F_a$ is a reticular $t$-${\cal P}$-${\cal K}$-stable unfolding of $F_a|_{t=0}$.
But the function germ $x_1^2+a x_1x_2+x_2^2\in {\mathfrak M}(2;0)^2$ has a modality.
This means that we can not use Theorem \ref{genericclass:tth} for the case $r=2$.
We require some equivalence relation of reticular Legendrian unfoldings 
which is weaker than the ${\cal P}$-Legendrian equivalence.\\

We give the figures of all generic bifurcations of wavefronts on a boundary.
Bifurcations of types ${}^0X$ do not occur.
We give figures of bifurcation of types ${}^1X$ at times $t<0,\ t=0,\ t>0$ 
respectively.
For example, the bifurcation of the type ${}^1B_3$ has
three bifurcations, that is one bifurcation in $2$D and two bifurcations in $3$D
(see Theorem \ref{classfunct:tth}).
\begin{figure}[htbp]
 \begin{minipage}{0.48\hsize}
  \begin{center}
    \includegraphics[width=5cm,height=5cm]{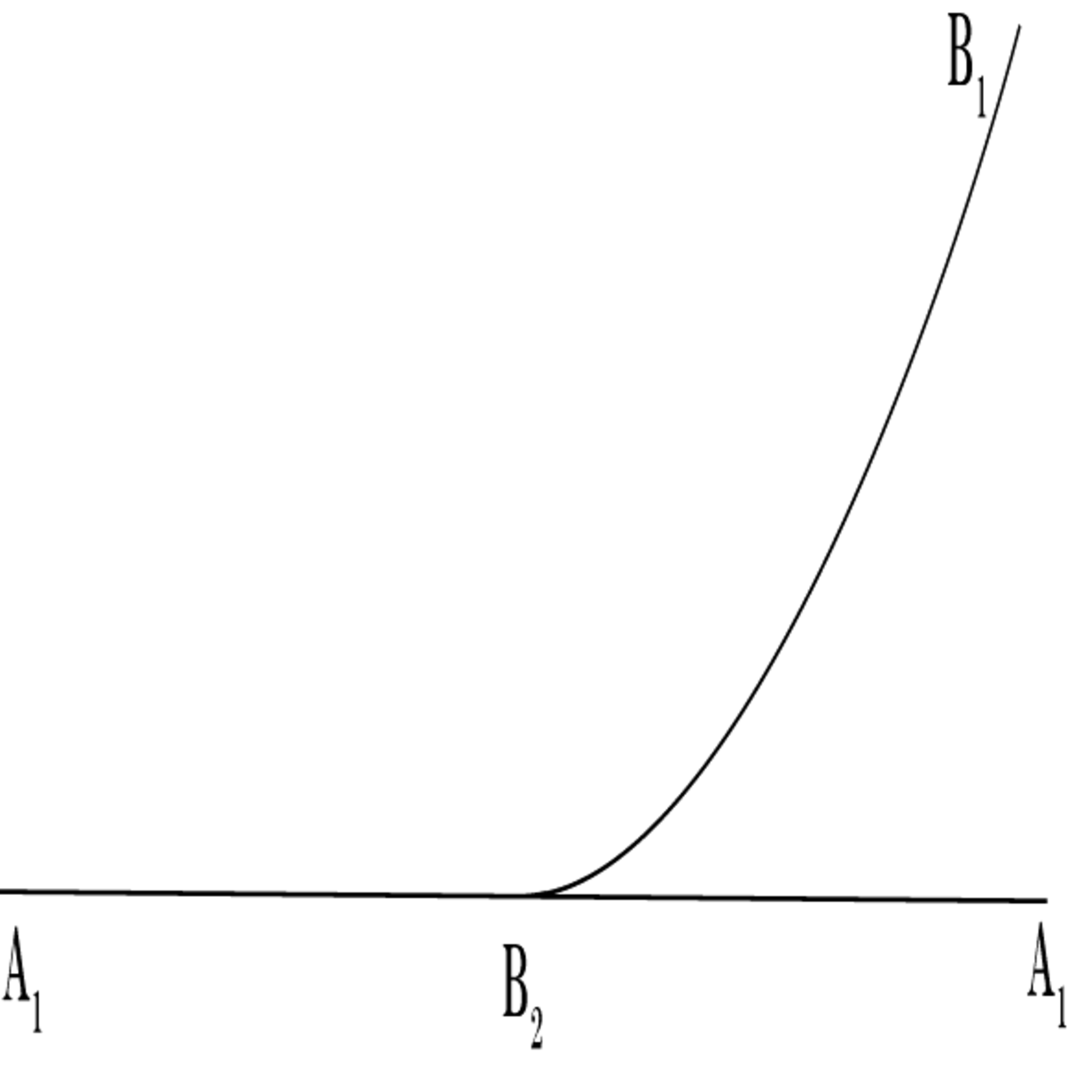}
  \end{center}
\caption{${}^0B_2$}
 \end{minipage}
 \begin{minipage}{0.48\hsize}
   \begin{center}
     \includegraphics[width=5cm,height=5cm]{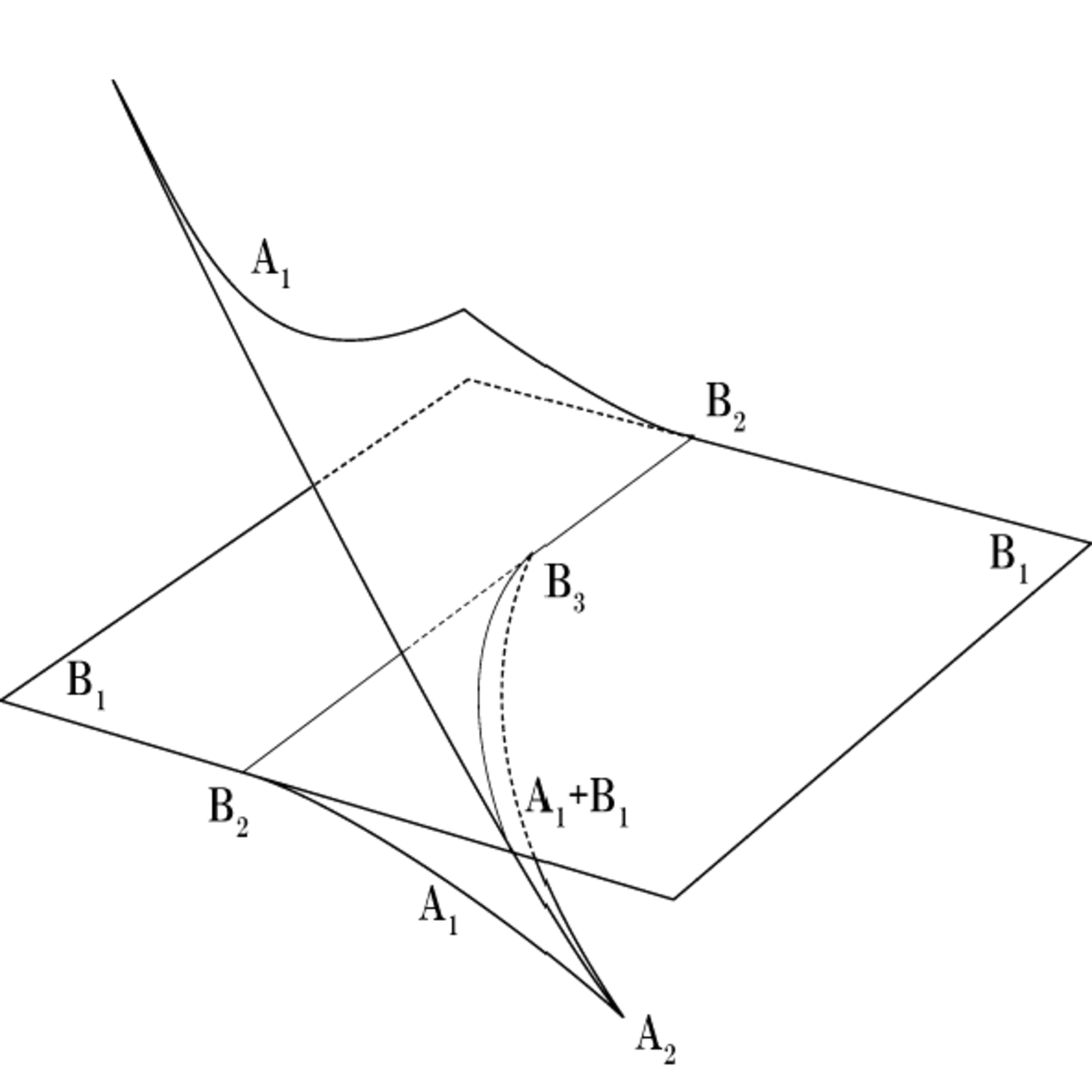}
   \end{center}
\caption{${}^0B_3$}
 \end{minipage}
\end{figure}

\begin{figure}[htbp]
 \begin{minipage}{0.48\hsize}
  \begin{center}
    \includegraphics[width=5cm,height=5cm]{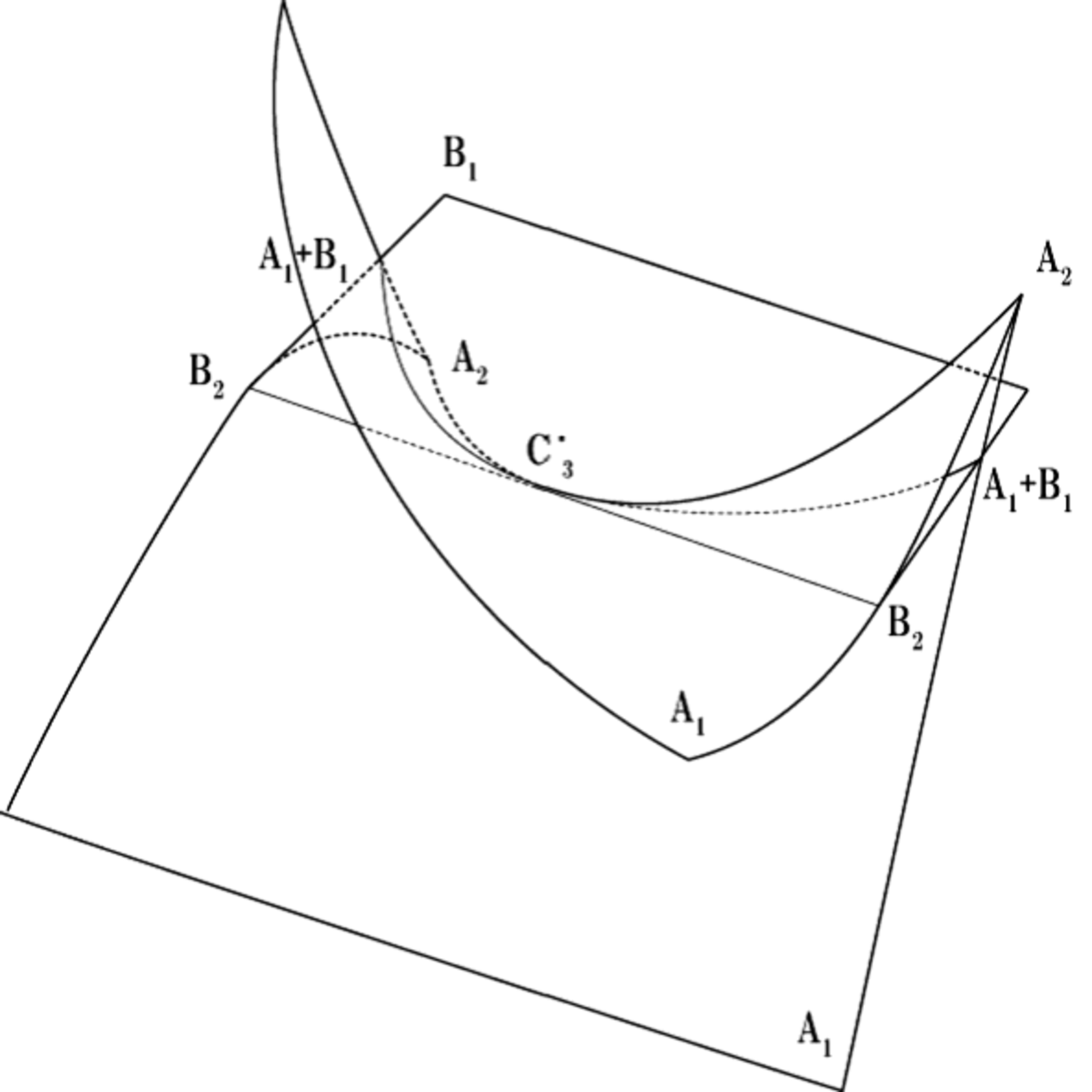}
  \end{center}
\caption{${}^0C^-_3$}
 \end{minipage}
 \begin{minipage}{0.48\hsize}
   \begin{center}
     \includegraphics[width=5cm,height=5cm]{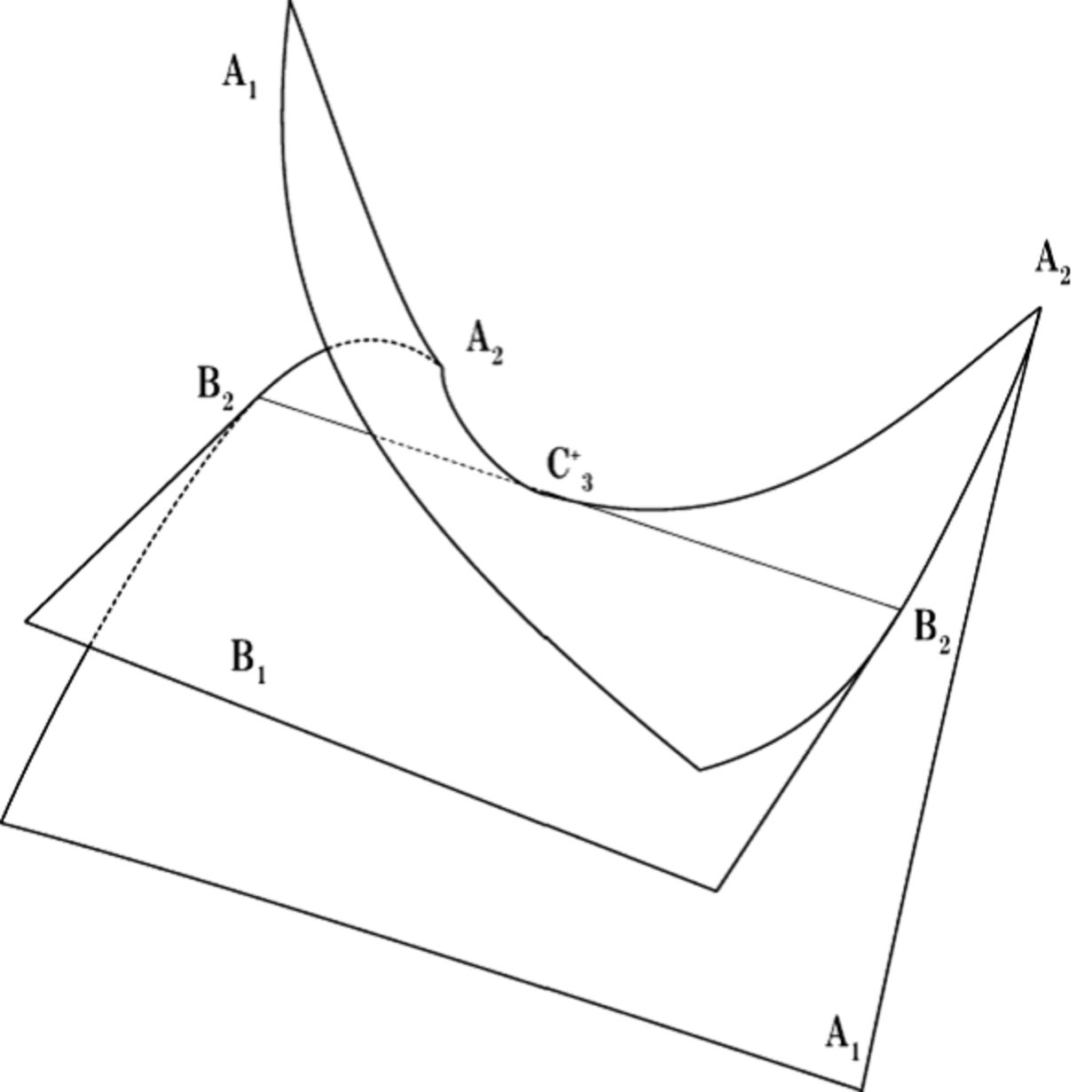}
   \end{center}
\caption{${}^0C^+_3$}
 \end{minipage}
\end{figure}

\begin{figure}[htbp]
 \begin{minipage}{0.30\hsize} 
  \begin{center}
    \includegraphics[width=3cm,height=3cm]{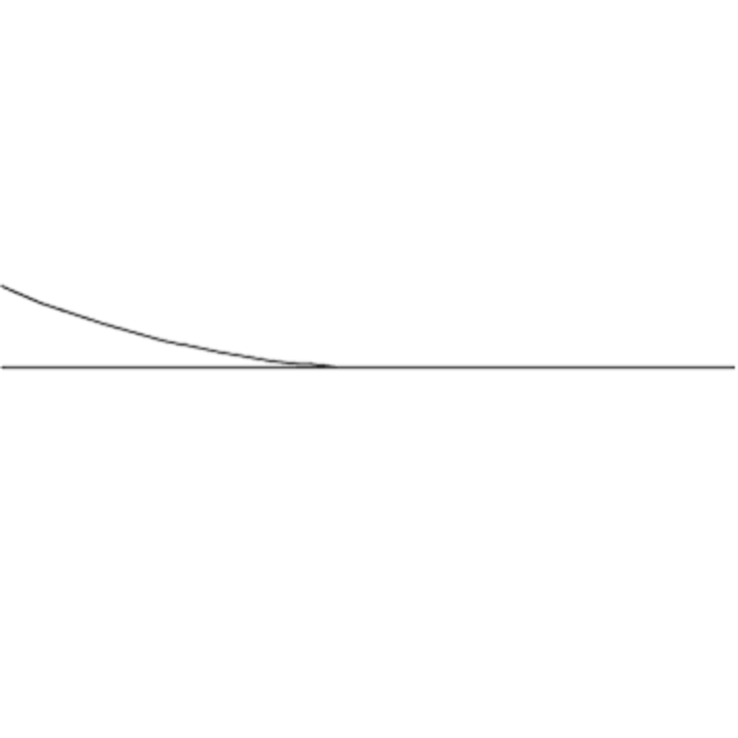}
  \end{center}
 \end{minipage}
$\leftrightarrow $
 \begin{minipage}{0.30\hsize}
   \begin{center}
     \includegraphics[width=3cm,height=3cm]{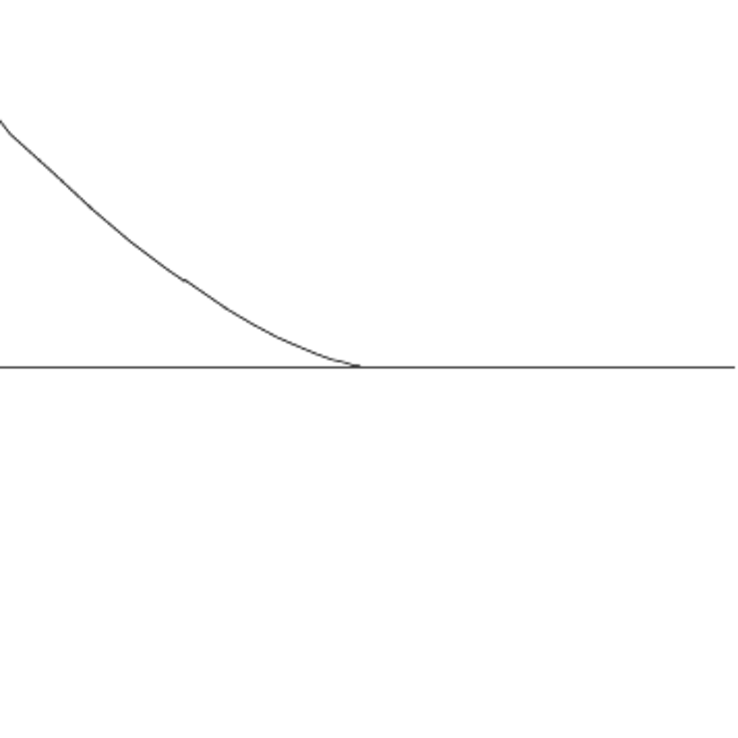}
   \end{center}
 \end{minipage}
$\leftrightarrow $
 \begin{minipage}{0.30\hsize}
  \begin{center}
  \includegraphics[width=3cm,height=3cm]{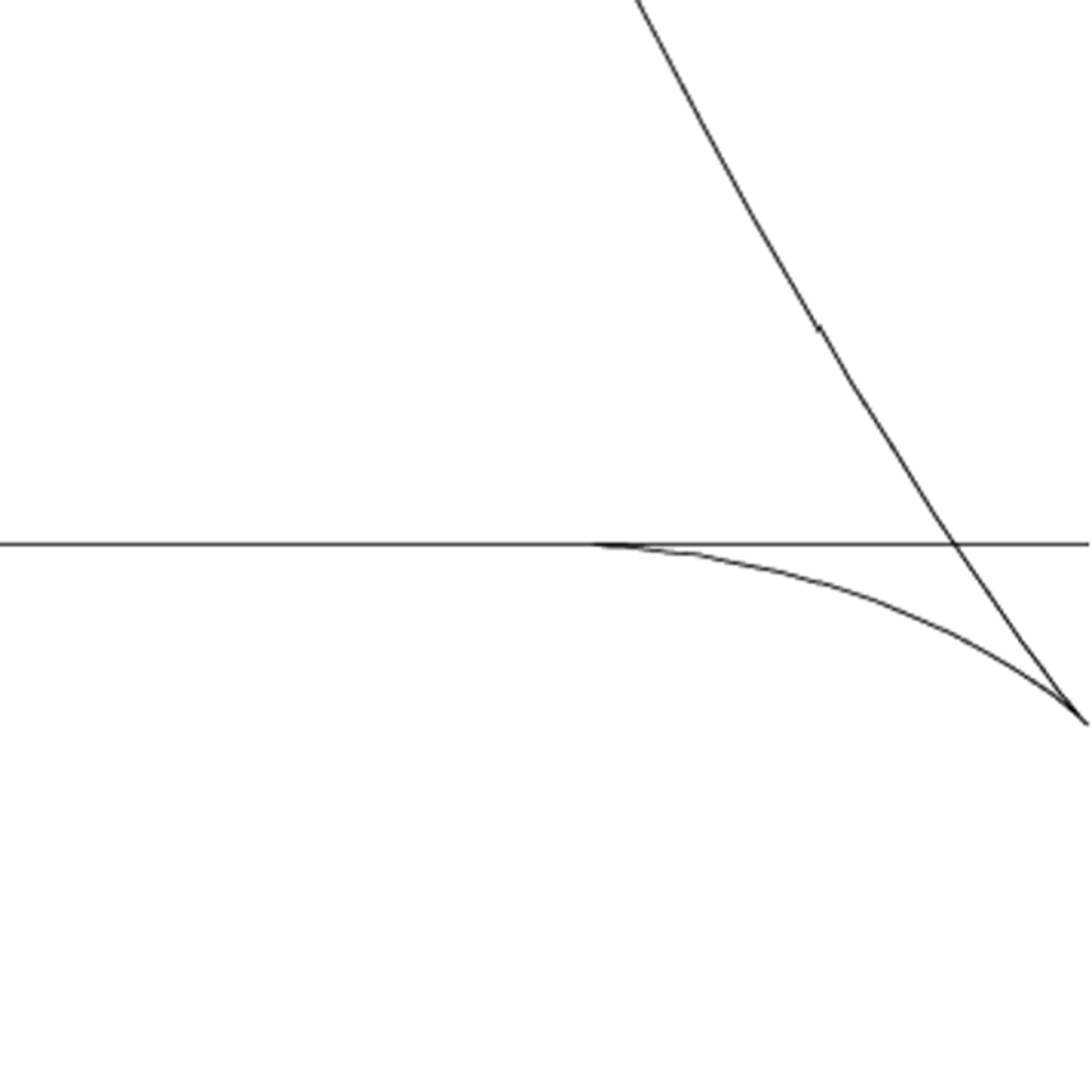}
  \end{center}
\end{minipage}\\
 \begin{minipage}{0.30\hsize}
  \begin{center}
   \includegraphics[width=3cm,height=3cm]{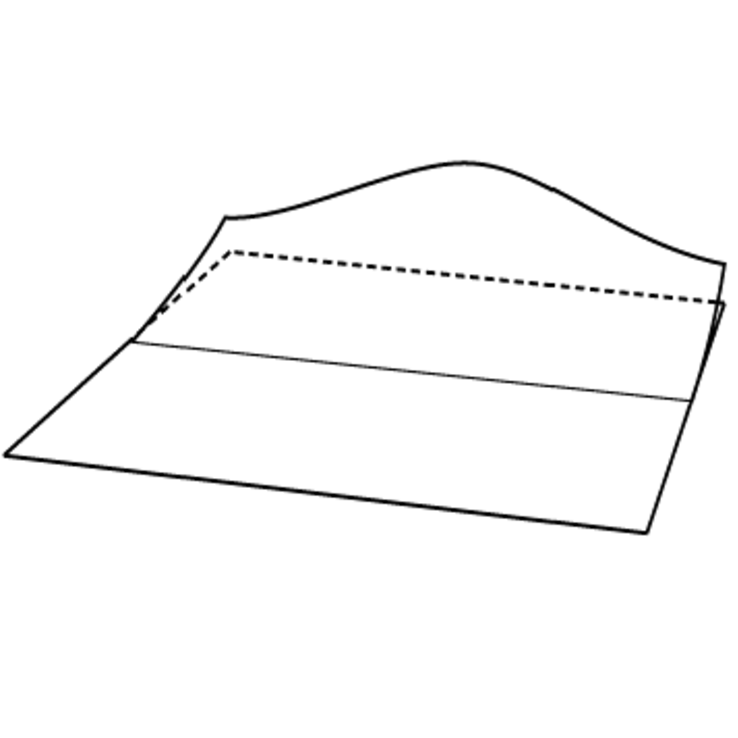}
  \end{center}
 \end{minipage}
$\leftrightarrow $
 \begin{minipage}{0.30\hsize}
   \begin{center}
     \includegraphics[width=3cm,height=3cm]{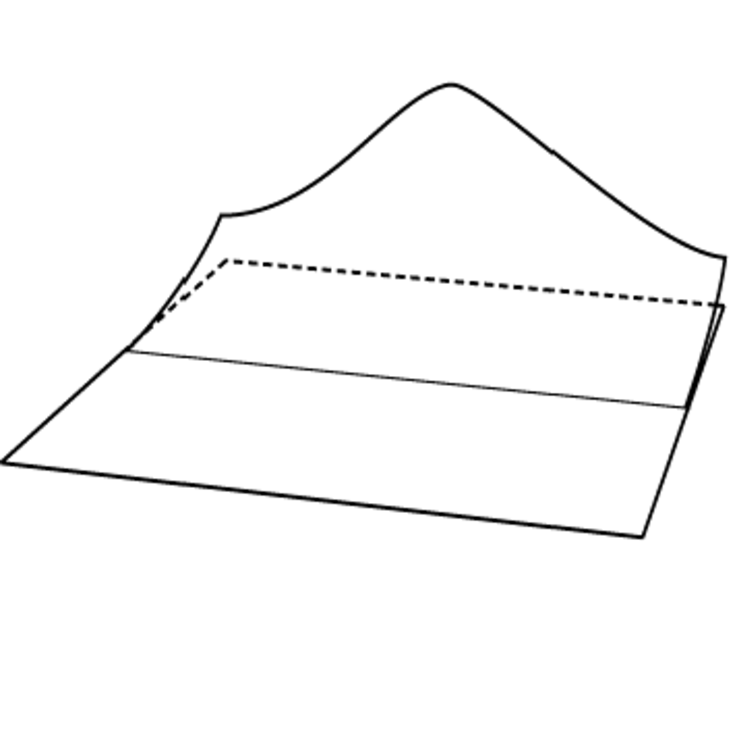}
   \end{center}
 \end{minipage}
$\leftrightarrow $
 \begin{minipage}{0.30\hsize}
  \begin{center}
  \includegraphics[width=3cm,height=3cm]{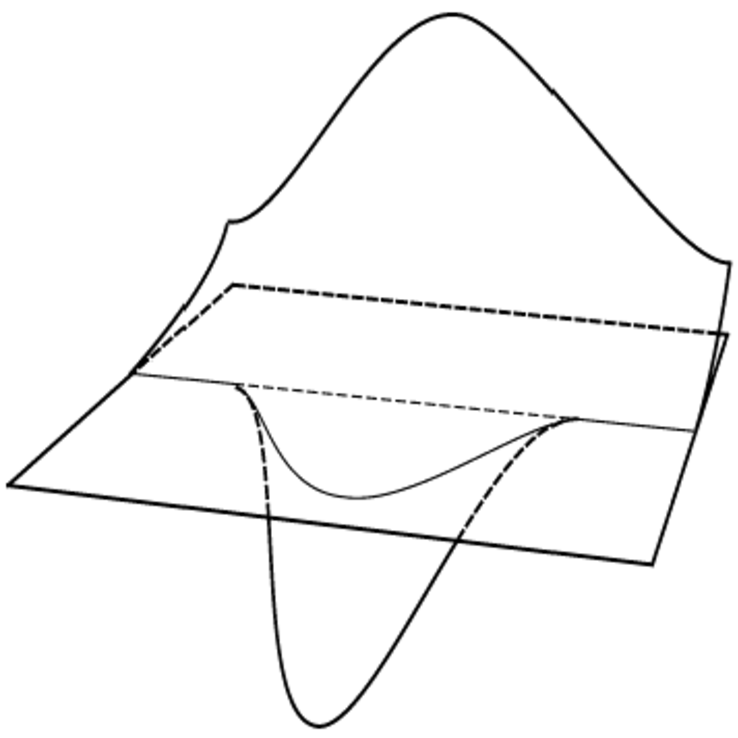}
  \end{center}
\end{minipage}
 \begin{minipage}{0.30\hsize}
  \begin{center}
    \includegraphics[width=3cm,height=3cm]{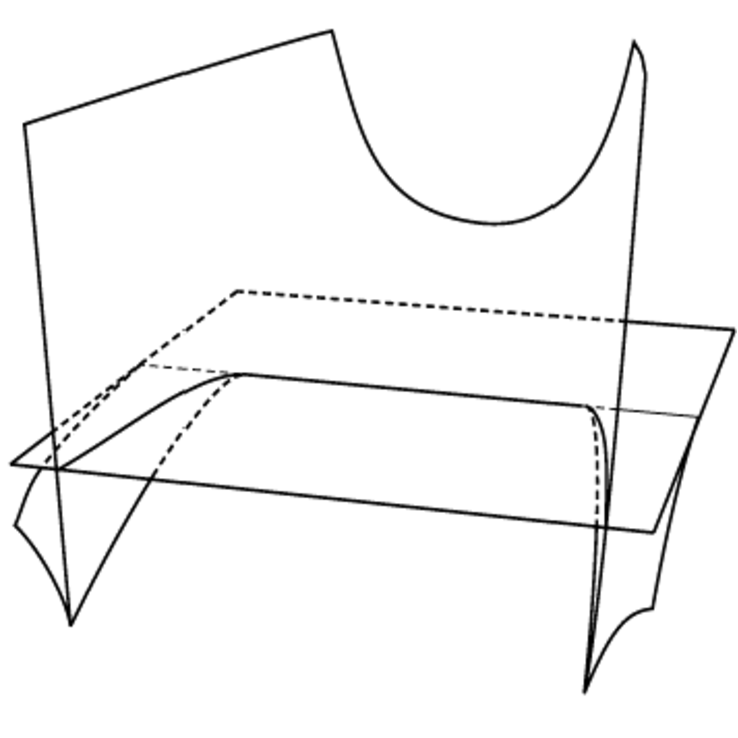}
  \end{center}
 \end{minipage}
$\leftrightarrow $
 \begin{minipage}{0.30\hsize}
   \begin{center}
     \includegraphics[width=3cm,height=3cm]{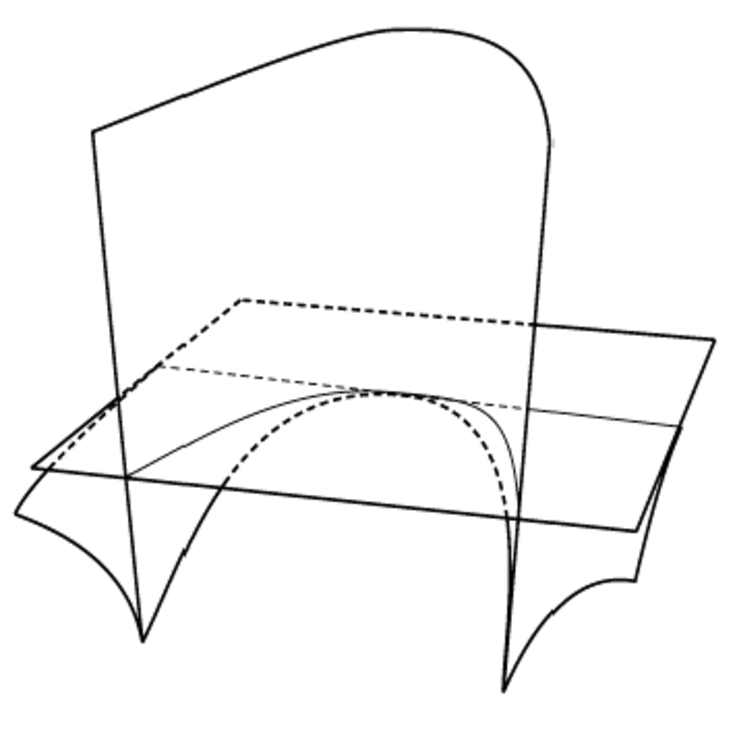}
   \end{center}
 \end{minipage}
$\leftrightarrow $
 \begin{minipage}{0.30\hsize}
  \begin{center}
  \includegraphics[width=3cm,height=3cm]{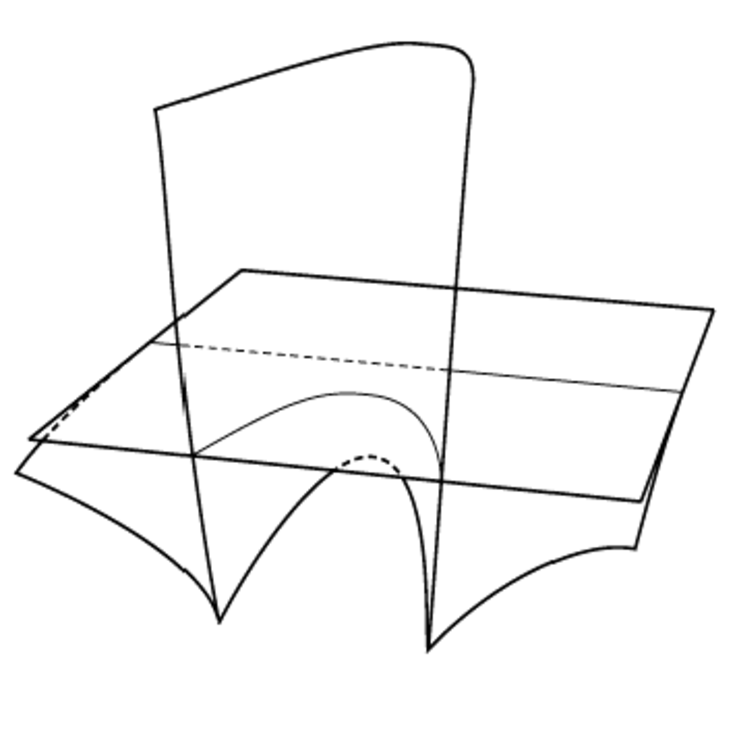}
  \end{center}
\end{minipage}
 \caption{${}^1B_3$}
\end{figure}

\begin{figure}[htbp]
 \begin{minipage}{0.30\hsize}
  \begin{center}
    \includegraphics[width=3cm,height=3cm]{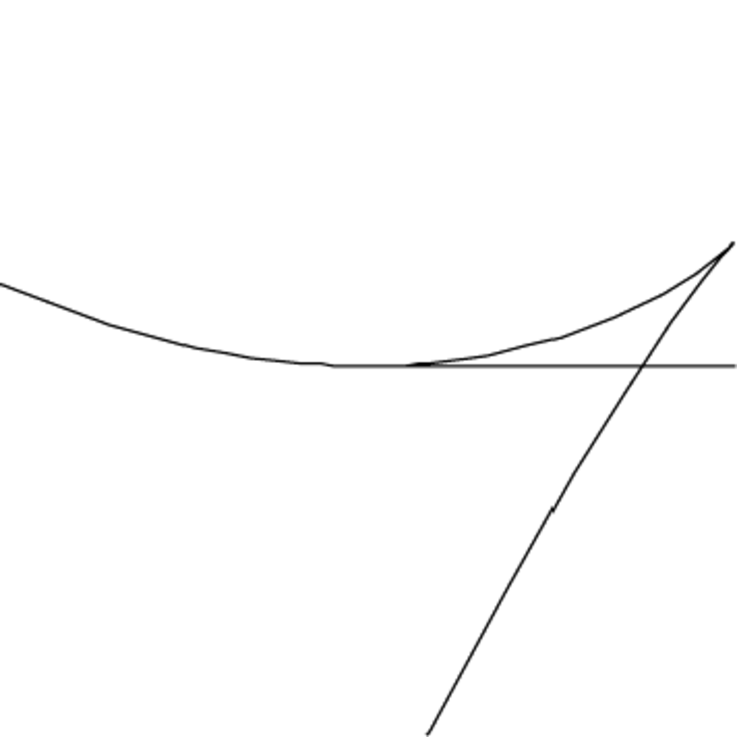}
  \end{center}
 \end{minipage}
$\leftrightarrow $
 \begin{minipage}{0.30\hsize}
   \begin{center}
     \includegraphics[width=3cm,height=3cm]{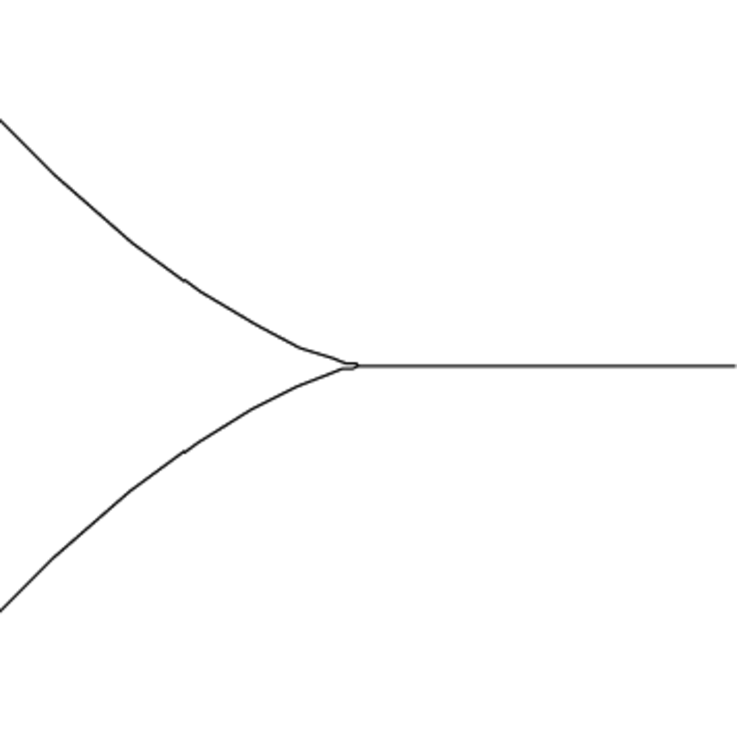}
   \end{center}
 \end{minipage}
$\leftrightarrow $
 \begin{minipage}{0.30\hsize}
  \begin{center}
  \includegraphics[width=3cm,height=3cm]{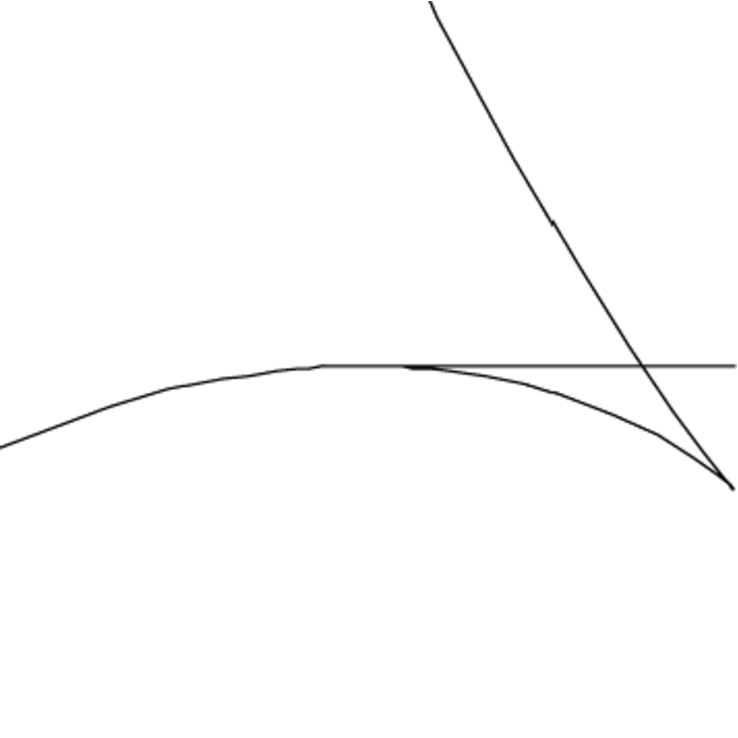}
  \end{center}
\end{minipage}\vspace{2mm}\\
 \begin{minipage}{0.30\hsize}
  \begin{center}
    \includegraphics[width=3cm,height=3cm]{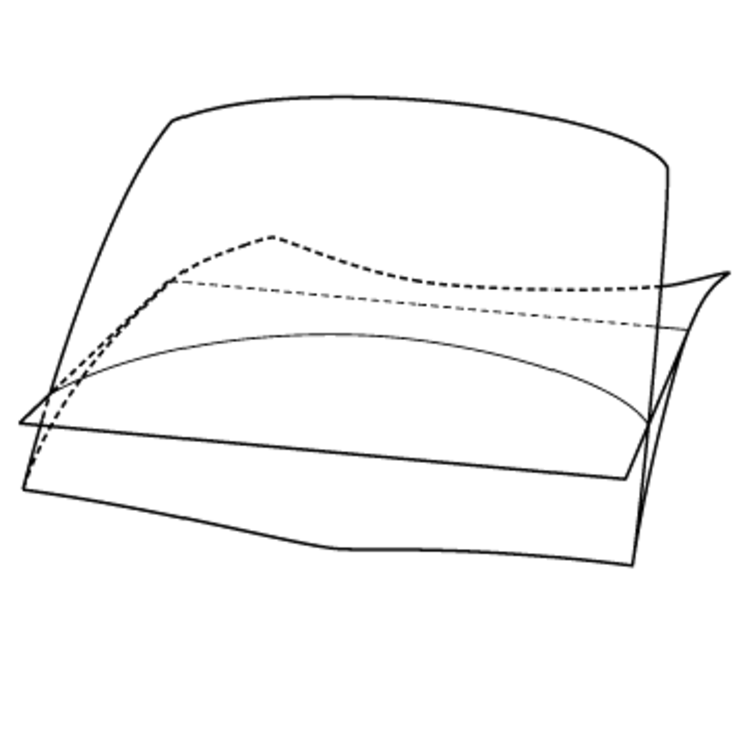}
  \end{center}
 \end{minipage}
$\leftrightarrow $
 \begin{minipage}{0.30\hsize}
   \begin{center}
     \includegraphics[width=3cm,height=3cm]{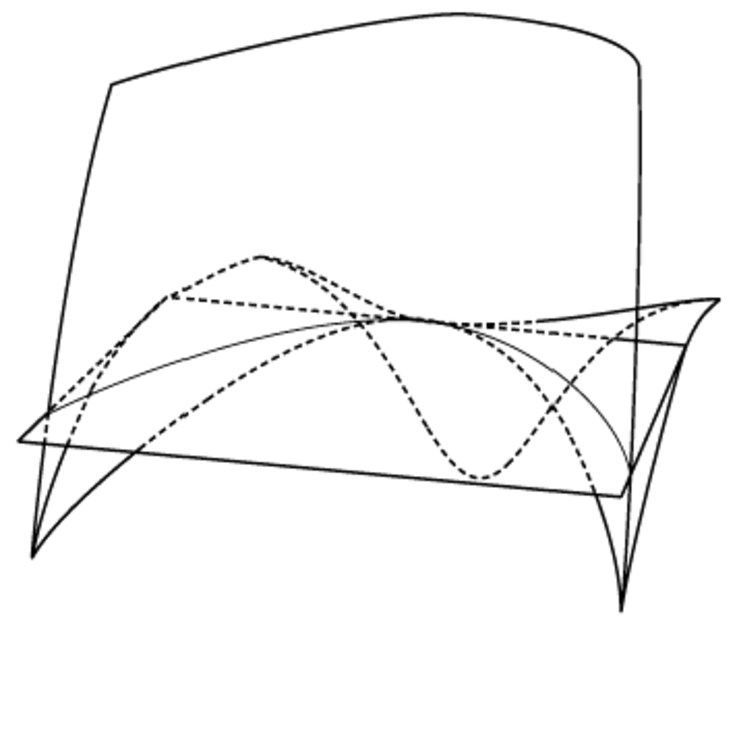}
   \end{center}
 \end{minipage}
$\leftrightarrow $
 \begin{minipage}{0.30\hsize}
  \begin{center}
  \includegraphics[width=3cm,height=3cm]{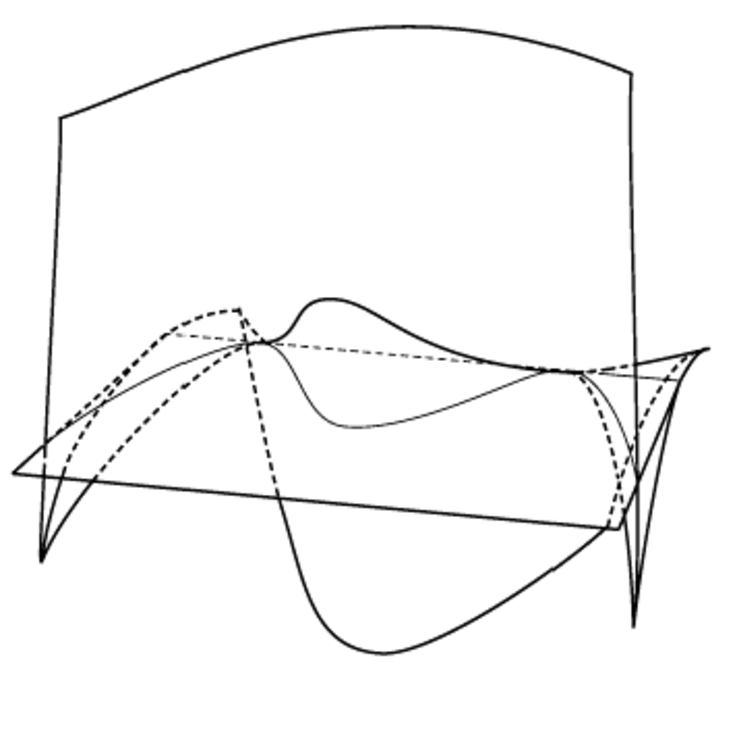}
  \end{center}
\end{minipage}
\caption{${}^1C^-_3$}
\end{figure}

\begin{figure}[htbp]
 \begin{minipage}{0.30\hsize}
  \begin{center}
    \includegraphics[width=3cm,height=3cm]{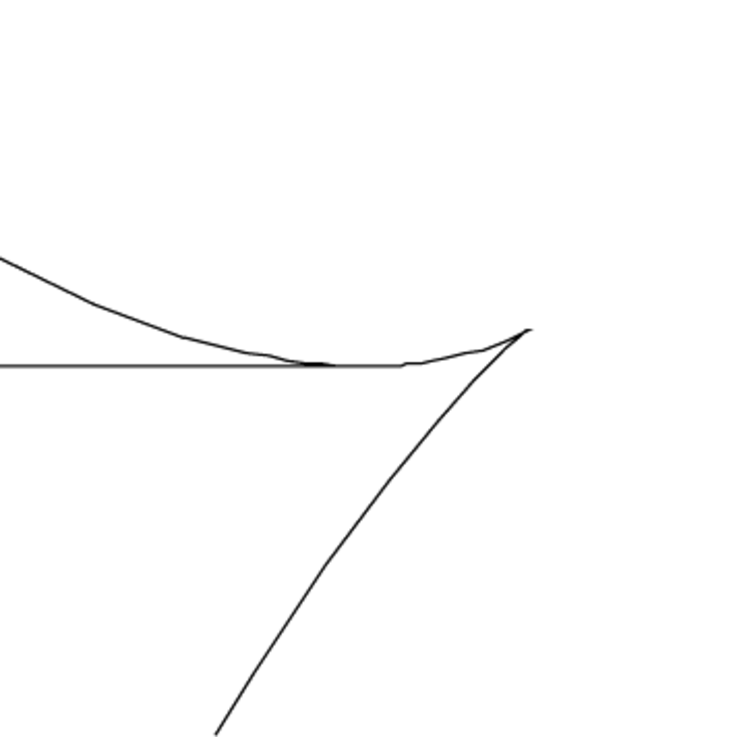}
  \end{center}
 \end{minipage}
$\leftrightarrow $
 \begin{minipage}{0.30\hsize}
   \begin{center}
     \includegraphics[width=3cm,height=3cm]{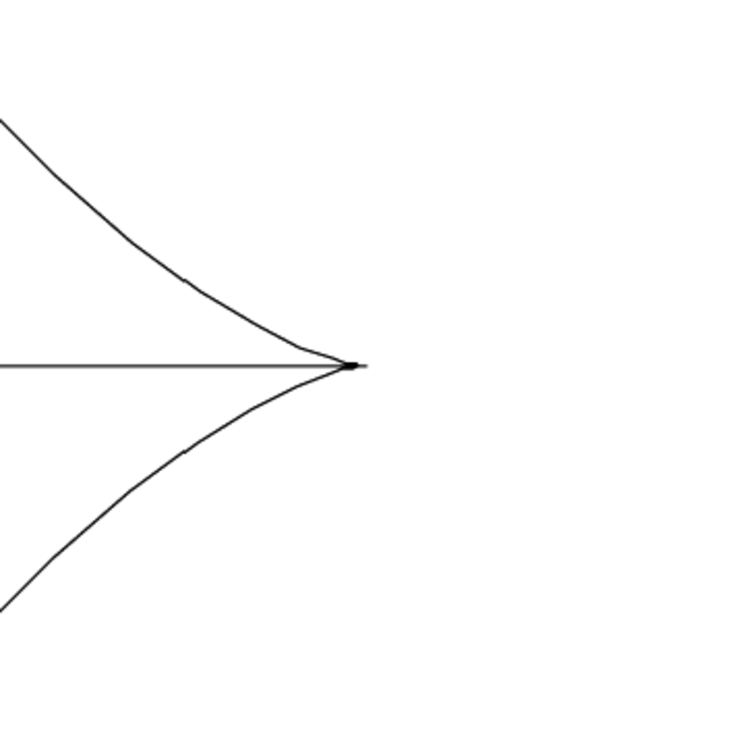}
   \end{center}
 \end{minipage}
$\leftrightarrow $
 \begin{minipage}{0.30\hsize}
  \begin{center}
  \includegraphics[width=3cm,height=3cm]{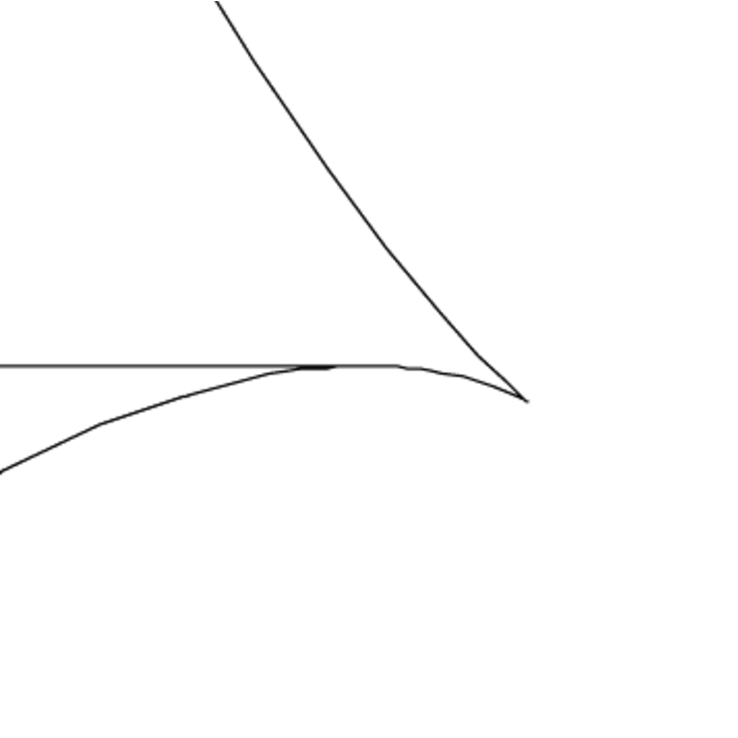}
  \end{center}
\end{minipage}\vspace{2mm} \\
 \begin{minipage}{0.30\hsize}
  \begin{center}
    \includegraphics[width=3cm,height=3cm]{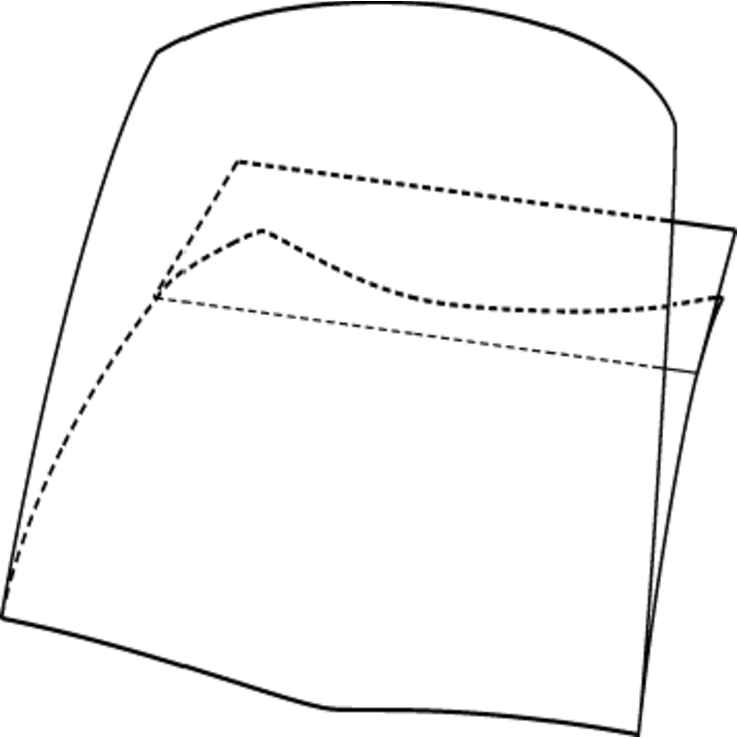}
  \end{center}
 \end{minipage}
$\leftrightarrow $
 \begin{minipage}{0.30\hsize}
   \begin{center}
     \includegraphics[width=3cm,height=3cm]{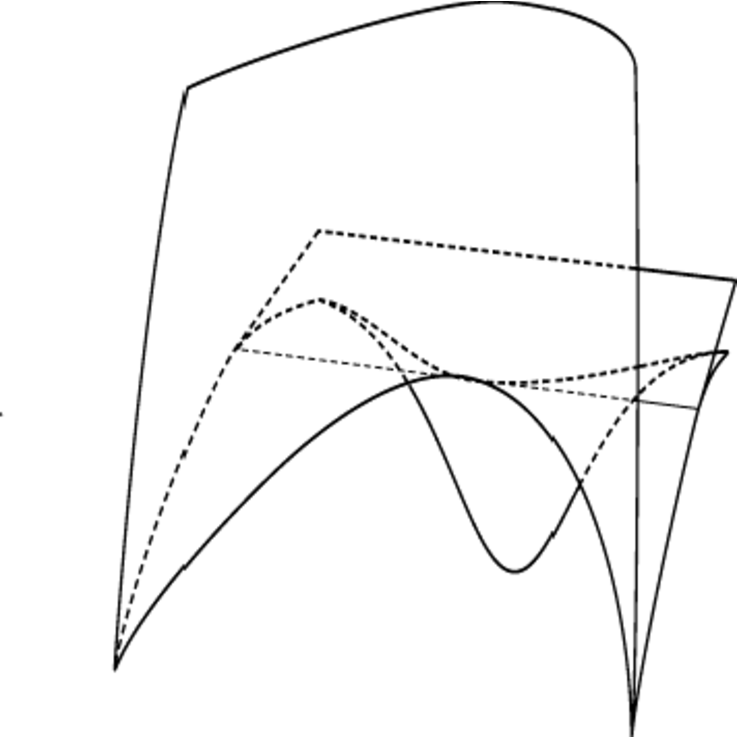}
   \end{center}
 \end{minipage}
$\leftrightarrow $
 \begin{minipage}{0.30\hsize}
  \begin{center}
  \includegraphics[width=3cm,height=3cm]{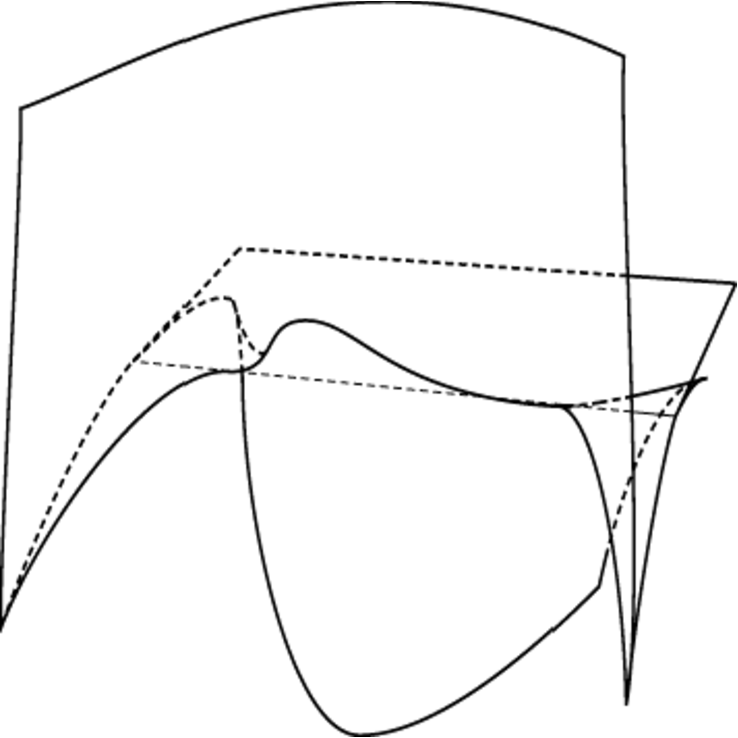}
  \end{center}
\end{minipage}
 \caption{${}^1C^+_3$}
\end{figure}
\begin{figure}[htbp]
 \begin{minipage}{0.30\hsize}
  \begin{center}
    \includegraphics[width=3cm,height=3cm]{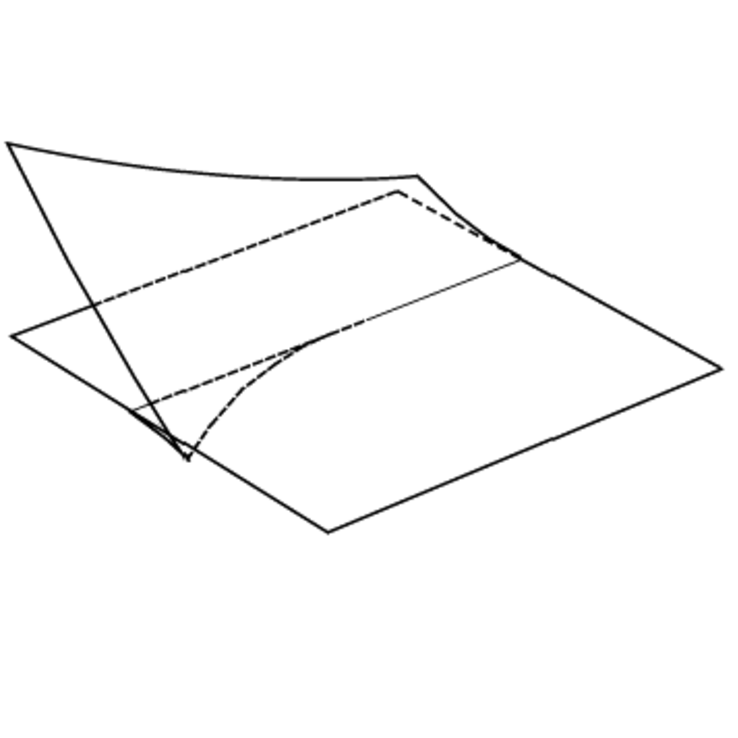}
  \end{center}
 \end{minipage}
$\leftrightarrow $
 \begin{minipage}{0.30\hsize}
   \begin{center}
     \includegraphics[width=3cm,height=3cm]{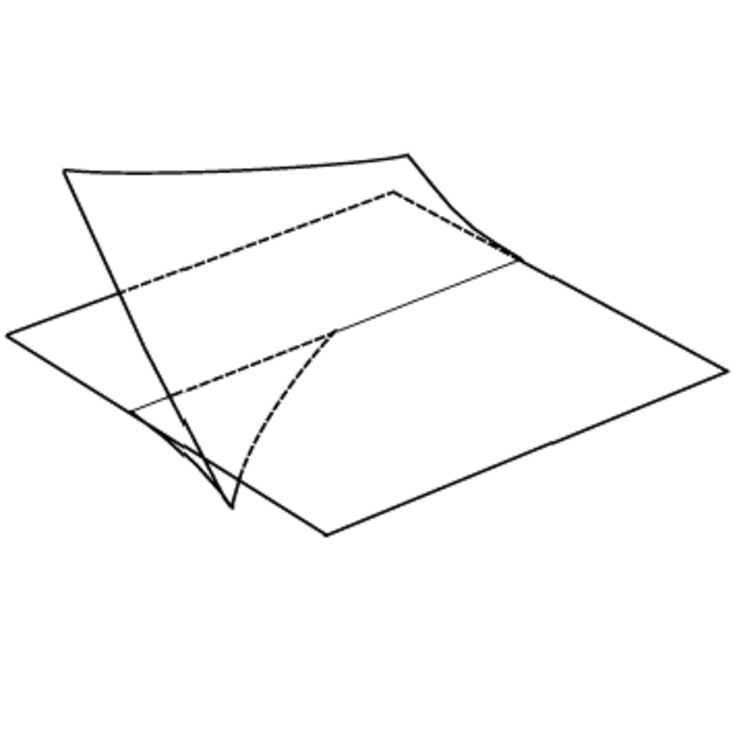}
   \end{center}
 \end{minipage}
$\leftrightarrow $
 \begin{minipage}{0.30\hsize}
  \begin{center}
  \includegraphics[width=3cm,height=3cm]{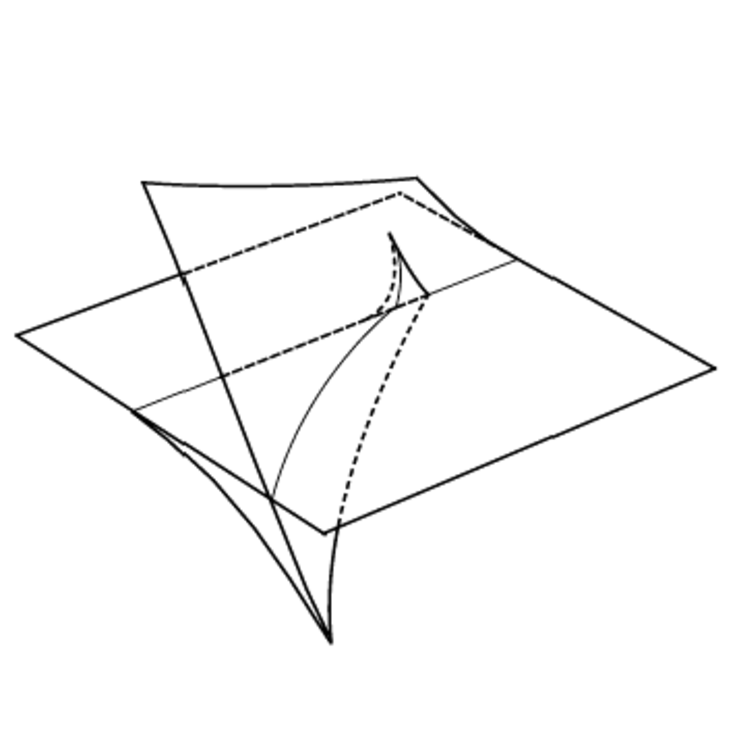}
  \end{center}
\end{minipage}
 \caption{${}^1B_4$}
\end{figure}
\begin{figure}[htbp]
 \begin{minipage}{0.30\hsize}
  \begin{center}
    \includegraphics[width=3cm,height=3cm]{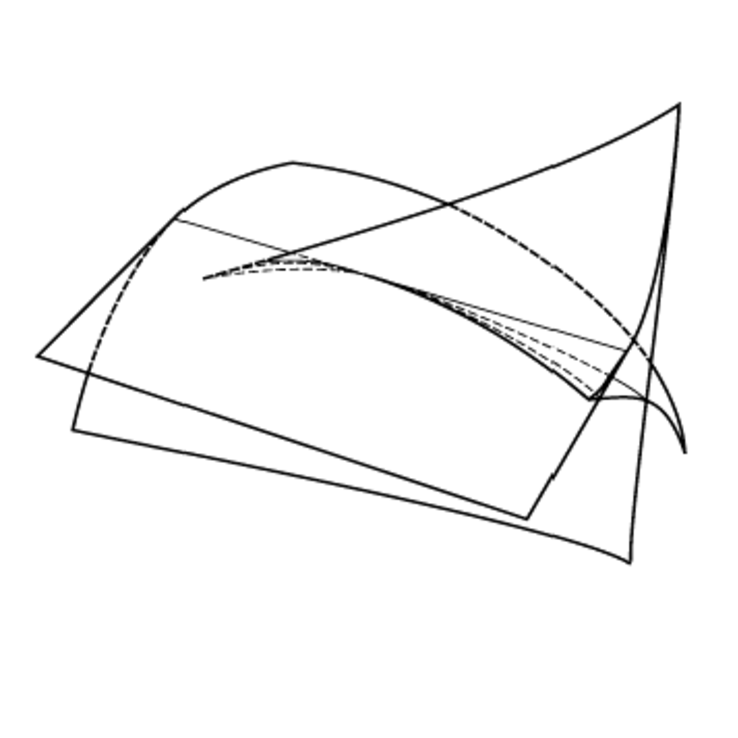}
  \end{center}
 \end{minipage}
$\leftrightarrow $
 \begin{minipage}{0.30\hsize}
   \begin{center}
     \includegraphics[width=3cm,height=3cm]{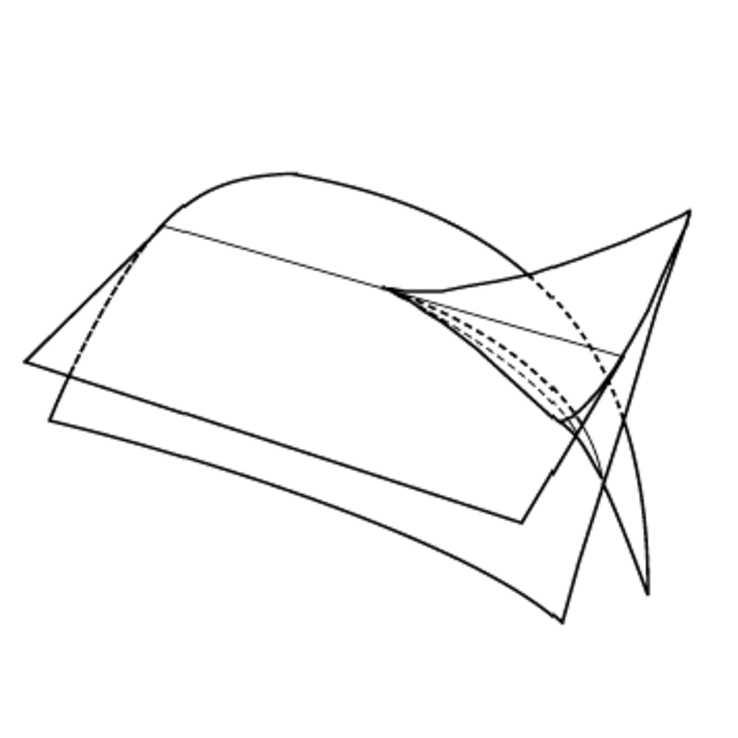}
   \end{center}
 \end{minipage}
$\leftrightarrow $
 \begin{minipage} {0.30\hsize}
  \begin{center}
  \includegraphics[width=3cm,height=3cm]{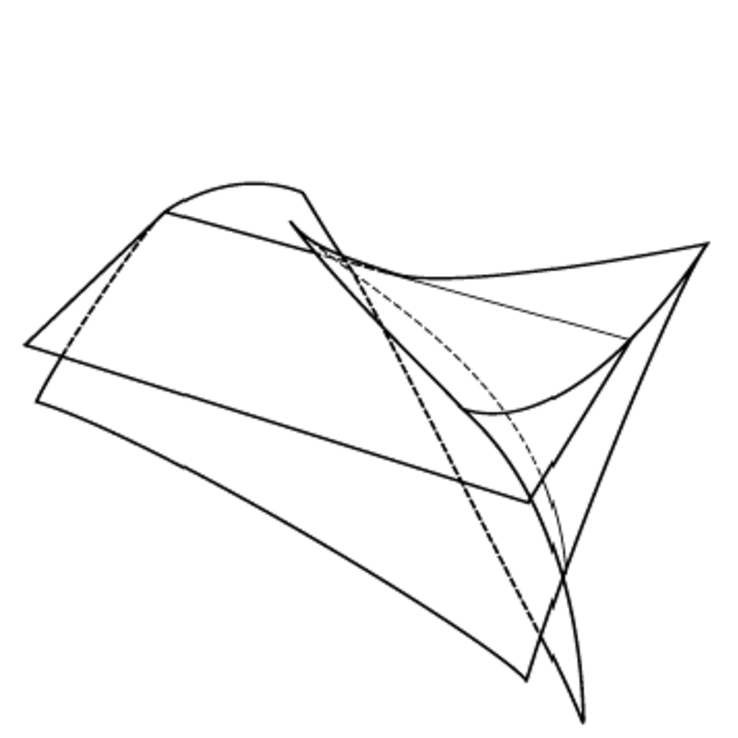}
  \end{center}
\end{minipage}
 \caption{${}^1C_4$}
\end{figure}
\begin{figure}[htbp]
 \begin{minipage}{0.30\hsize}
  \begin{center}
    \includegraphics[width=3cm,height=3cm]{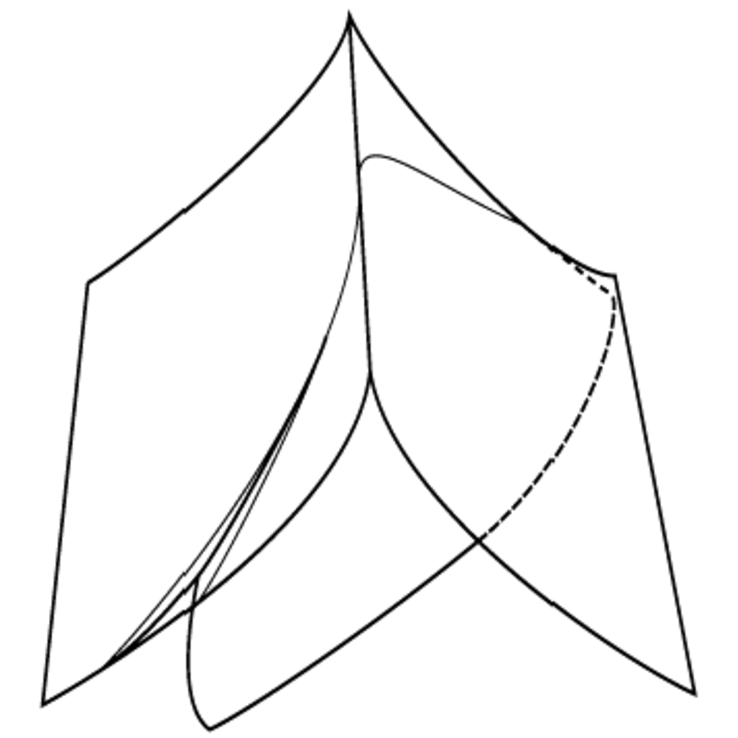}
  \end{center}
 \end{minipage}
$\leftrightarrow $
 \begin{minipage}{0.30\hsize}
   \begin{center}
     \includegraphics[width=3cm,height=3cm]{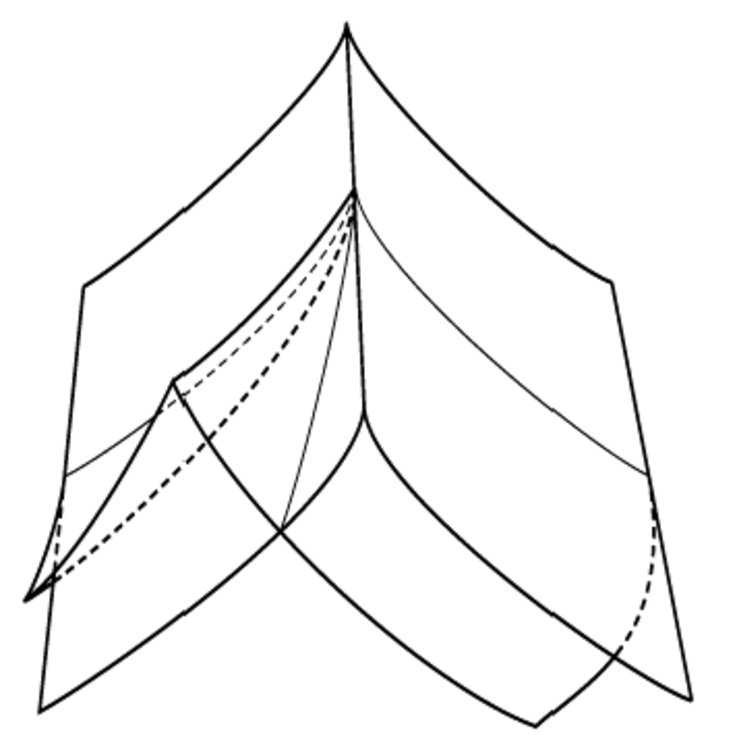}
   \end{center}
 \end{minipage}
$\leftrightarrow $
 \begin{minipage} {0.30\hsize}
  \begin{center}
  \includegraphics[width=3cm,height=3cm]{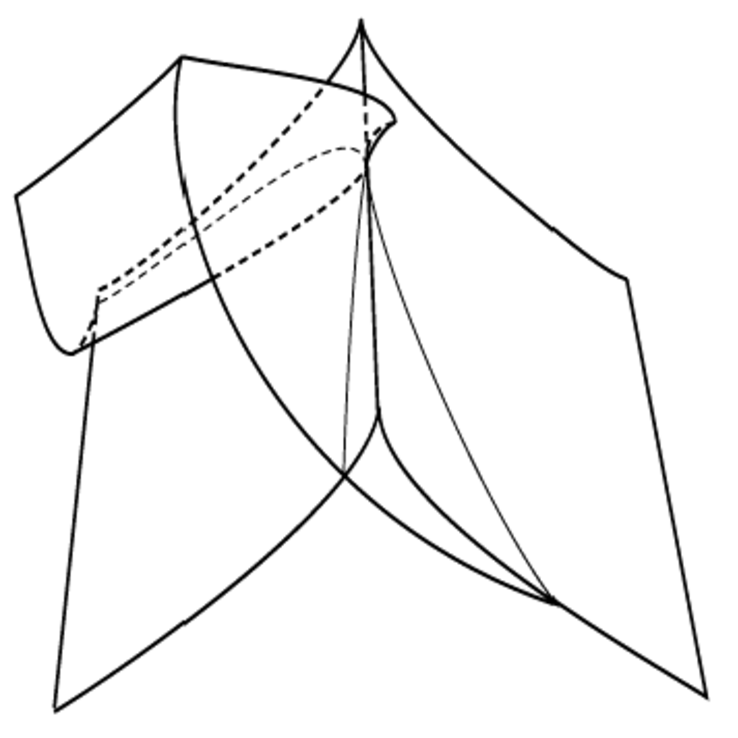}
  \end{center}
\end{minipage}
 \caption{${}^1F_4$}
\end{figure}          

\bibliographystyle{plain}
\bibliography{Bibt}
\end{document}